\documentclass[a4paper,10pt]{article}
\usepackage[english]{babel}
\usepackage[utf8]{inputenc}

\usepackage{amsmath}
\usepackage{amsfonts}
\usepackage{amssymb}
\usepackage{amsthm}
\usepackage{graphicx}
\usepackage{pdfpages}

\newcommand{\mytext}[1]{ \: \textrm{#1} \: }
\newcommand{\mysetdescr}[2]{\left\{ { #1 \: \left| \: #2 \right. } \right\} }
\newcommand{\mysetdef}[3]{#1 \; \equiv \;\mysetdescr{#2}{#3} }
\newcommand{\myfr}[1]{ \mathfrak{#1} }
\newcommand{\mydownarrow}{{\downarrow \,}}
\newcommand{\myuparrow}{{\uparrow \,}}
\newcommand{\myouparrow}{\uparrow_{_{_{\!\!\!\circ}}}}
\newcommand{\myodownarrow}{\downarrow^{^{\!\!\!\circ}}}

\newcommand{\myN}{\mathbb{N}}
\newcommand{\myNk}[1]{\underline {#1}}
\newcommand\mytimes{{\times}}

\newcommand\myurbild[1]{\overline{#1}^1}

\def\A{{\cal A}}

\def\E{{\cal E}}

\def\G{{\cal G}}
\def\H{{\cal H}}
\def\I{{\cal I}}
\def\J{{\cal J}}
\def\R{{\cal R}}
\def\S{{\cal S}}
\def\T{{\cal T}}

\newcommand{\mfP}{ \myfr{P} }
\newcommand{\mfE}{ \myfr{E} }

\newcommand{\myGxi}[1]{\G(#1)}
\newcommand{\Gxi}{\myGxi{\xi}}
\newcommand{\Gzet}{\myGxi{\zeta}}

\newcommand{\mygxix}[2]{G_{#1}(#2)}
\newcommand{\gxix}{\mygxix{\xi}{x}}
\newcommand{\grxix}{\mygxix{\rho(\xi)}{x}}

\newcommand{\gzex}{\mygxix{\zeta}{x}}

\newcommand{\myaxix}[2]{\alpha_{#1}(#2)}

\newcommand{\atsX}{\alpha_{\tau \circ \sigma}(x)}

\newcommand{\atSX}{\alpha_{\tau}(\sigma(x))}
\newcommand{\atSY}{\alpha_{\tau}(\sigma(y))}

\newcommand{\axix}{\myaxix{\xi}{x}}

\newcommand{\myapxix}[3]{\alpha_{#1, #2}(#3)}

\newcommand{\myapepxix}[3]{\alpha_{#1, \eta_{#1}(#2)}(#3)}

\newcommand{\apepxi}{\alpha_{P,\eta_P(\xi)}}
\newcommand{\myepxix}[3]{\eta_{#1}(#2)(#3)}

\newcommand{\mf}[1]{\mathfrak{ #1 }}
\newcommand{\fa}{\mf{a}}
\newcommand{\fb}{\mf{b}}
\newcommand{\fc}{\mf{c}}

\newcommand{\eps}{\epsilon}

\newcommand{\arie}{P \in \mfP_r$, $\xi \in \H(P,R)$, $x \in P}
\newcommand{\ariekl}{\xi \in \H(P,R)$, $x \in P}

\newcommand{\subseteqA}{\, \subseteq \,}
\newcommand{\inA}{\in \,}

\def\BP{\begin{proof}}
\def\EP{\end{proof}}

\DeclareMathOperator{\Aut}{Aut}
\DeclareMathOperator{\Emb}{Emb}

\begin{document}

\theoremstyle{plain}
\newtheorem{condition}{Condition}
\newtheorem{theorem}{Theorem}
\newtheorem{definition}{Definition}
\newtheorem{corollary}{Corollary}
\newtheorem{lemma}{Lemma}
\newtheorem{proposition}{Proposition}

\title{\bf Strong G-schemes and strict homomorphisms}
\author{\sc Frank a Campo}
\date{\small Seilerwall 33, D 41747 Viersen, Germany\\
{\sf acampo.frank@gmail.com}}



\maketitle

\begin{abstract}
\noindent Let $\mfP_r$ be a representation system of the non-isomorphic finite posets, and let $\H(P,Q)$ be the set of order homomorphisms from $P$ to $Q$. For finite posets $R$ and $S$, we write $R \sqsubseteq_G S$ iff, for every $P \in \mfP_r$, a one-to-one mapping $\rho_P : \H(P,R) \rightarrow \H(P,S)$ exists which fulfills a certain regularity condition. It is shown that $R \sqsubseteq_G S$ is equivalent to $\# \S(P,R) \leq \# \S(P,S)$ for every finite posets $P$, where $\S(P,Q)$ is the set of strict order homomorphisms from $P$ to $Q$. In consequence, $\# \S(P,R) = \# \S(P,S)$ holds for every finite posets $P$ iff $R$ and $S$ are isomorphic. A sufficient condition is derived for $R \sqsubseteq_G S$ which needs the inspection of a finite number of posets only. Additionally, a method is developed which facilitates for posets $P + Q$ (direct sum) the construction of posets $T$ with $P + Q \sqsubseteq_G A + T$, where $A$ is a convex subposet of $P$.
\newline

\end{abstract}

\section{Introduction} \label{sec_introduction}

Garret Birkhoff \cite{Birkhoff_1937,Birkhoff_1942} published two articles in 1937 and 1942 in which he introduced for posets $P$ and $Q$ their direct sum $P + Q$, their product $P \times Q$, and the homomorphism set $\H(P,Q)$, together with their (later on) usual partial order relations. His purpose was to unify ordinal and cardinal arithmetic, but doing so, he opened the rich field of ``order arithmetic''. For surveys, see J\'{o}nsson \cite{Jonsson_1982_1} from 1982, Duffus \cite{Duffus_1984} from 1984, and McKenzie \cite{McKenzie_2003} from 2003.

The operations introduced by Birkhoff were extended by Day \cite{Day_1945} to more general relations in 1945. As an important step, Hashimoto \cite{Hashimoto_1948,Hashimoto_1951} proved in 1948 and 1951, that two product representations of a connected poset have always a common refinement. Lov\'{a}sz \cite{Lovasz_1967} showed in 1967, that there is an in-depth connection between the structure of posets and the cardinalities of homomorphism sets related to them. He applied the result in showing a cancellation rule for the product of finite posets: $P \times R \simeq P \times S \Rightarrow R \simeq S$, and extended in \cite{Lovasz_1971} his result to more general structures in 1971. (Here and for the rest of the overview, all posets are non-empty.)

Birkhoff \cite{Birkhoff_1942} conjectured already in 1942, that $\H(P,R) \simeq \H(P,S)$ implies $R \simeq S$ for finite posets $P, R$, and $S$. This so-called ``cancellation problem'', moved into focus at the end of the seventies. Progress was made for many cases. For a chain $P$, the cancellation rule was proven by Bergman et al.\ \cite{Bergman_etal_1977} in 1977. In 1978, Duffus et al.\ \cite{Duffus_etal_1978} worked with $R$ being a lattice, and in 1980, Wille \cite{Wille_1980} contributed to this topic as well. For finite connected posets $P$ and $Q$, Duffus and Wille \cite{Duffus_Wille_1979} showed in 1979 that $\H(P,P) \simeq \H(Q,Q)$ implies $P \simeq Q$.

Duffus and Rival \cite{Duffus_Rival_1978} proved in 1978 a ``logarithmic property'', which transfers certain cancellation problems for homomorphism sets to cancellation problems for products of posets. Duffus \cite{Duffus_1984} developed this approach further in 1984.

In 1982, J\'{o}nsson \cite{Jonsson_1982_1} published an overview over the arithmetic of ordered sets, providing a rich spectrum of results and aspects. In the same year, J\'{o}nsson  \cite{Jonsson_1982_2} and J\'{o}nsson and McKenzie \cite{Jonsson_McKenzie_1982} presented several important results about order arithmetic and cancellation rules.

The starting signal for finally solving the problem was given by Farley \cite{Farley_1996} in 1996, with an article about the structure of the automorphism group of certain sets of homomorphisms. In 1999 and 2000, McKenzie \cite{McKenzie_1999,McKenzie_2000} proved that indeed $\H(P,R) \simeq \H(P,S)$ implies $R \simeq S$ for finite poests $P, R, S$. In 2003, McKenzie \cite{McKenzie_2003} published an additional paper about this subject, using a different approach. After this publication, the interest in homomorphism sets as such expired.

The present paper belongs to a series of articles \cite{aCampo_toappear_1,aCampo_toappear_2} which deal with homomorphism sets under the following aspect: {\em What is it in the structure of finite posets $R$ and $S$ that results in $\# \H(P,R) \leq \# \H(P,S)$ for every finite poset $P$?}. An early contribution to this topic is provided by Lov\'{a}sz \cite{Lovasz_1967} who showed in 1967, that $\# \H(P,R) = \# \H(P,S)$ for every finite poset $P$ is equivalent to $R \simeq S$.

Let $\mfP$ denote the class of all finite posets, and let $\mfP_r$ denote a representation system of the non-isomorphic posets in $\mfP$. In a previous article \cite{aCampo_toappear_1}, the concepts of the Hom-scheme and the G-scheme have been introduced: given finite posets $R$ and $S$, a {\em Hom-scheme} $\rho$ from $R$ to $S$ assigns a mapping $\rho_P : \H(P,R) \rightarrow \H(P,S)$ to every poset $P \in \mfP_r$. A Hom-scheme $\rho$ is called a {\em G-scheme} iff, for every poset $P \in \mfP_r$ and every $\ariekl$, the {\em connectivity component $K$ of the pre-image of $\xi(x)$ with $x \in K$} is identical with the connectivity component $K'$ of the pre-image of $\rho_P(\xi)(x)$ with $x \in K'$. A G-scheme is thus a Hom-scheme which obeys a regularity condition in the way how it maps homomorphisms.

A Hom-scheme / G-scheme is called {\em strong} iff $\rho_P$ is one-to-one for every $P \in \mfP_r$. The existence of a strong Hom-scheme from $R$ to $S$ is thus equivalent to $\# \H(P,R) \leq \# \H(P,S)$ for every $P \in \mfP$. However, the existence of a strong G-scheme has two aspects: the first one is $\# \H(P,R) \leq \# \H(P,S)$ for every poset $P \in \mfP$, the second one is the regularity condition.

In this paper it is shown, that also the strong G-schemes can be characterized by a single relation between cardinalities of homomorphism sets: in Theorem \ref{theo_GschemeOnStrict} in Section \ref{subsec_GschemeOnStrict}, we prove that the existence of a strong G-scheme from $R$ to $S$ is equivalent to $\# \S(P,R) \leq \# \S(P,S)$ for all $P \in \mfP$, where $\S(P,Q)$ is the set of {\em strict} homomorphisms from $P$ to $Q$. In consequence, $\# \S(P,R) = \# \S(P,S)$ holds for all finite posets $P$ iff $R$ and $S$ are isomorphic (Corollary \ref{coro_Lovaszrefinement}).

We write $R \sqsubseteq_G S$ iff a strong G-scheme exists from $R$ to $S$. If we want to apply Theorem \ref{theo_GschemeOnStrict} in the proof of $R \sqsubseteq_G S$ for given posets $R$ and $S$, we still have to work with the infinite number of posets $P$ contained in $\mfP_r$. Building on Theorem \ref{theo_GschemeOnStrict}, we achieve in Theorem \ref{theo_Transport} in Section \ref{sec_suffCond} a sufficient condition for $R \sqsubseteq_G S$ which needs the inspection of a {\em finite} number of {\em connected} posets only. We apply it in two ways. In Section \ref{subsec_proving_RGS}, we demonstrate on the basis of two examples that this theorem (if applicable) simplifies the proof of $R \sqsubseteq_G S$ considerably. Additionally, we present in Theorem \ref{theo_PQ_AT} in Section \ref{subsec_constr_Gschemes} a method which facilitates for posets of the form $P + Q$ (direct sum) the construction of posets $T$ with $P + Q \sqsubseteq_G A + T$, where $A$ is a convex subposet of $P$. Moreover, if $A$ is an antichain, there exists even a strong I-scheme from $P + Q$ to $A + T$. (A strong I-scheme is a strong G-scheme fulfilling an additional regularity condition; strong I-schemes are treated in \cite{aCampo_toappear_1}.)

As mentioned above, this paper belongs to a series of articles \cite{aCampo_toappear_1,aCampo_toappear_2} dealing with the question, what it is in the structure of finite posets $R$ and $S$ that results in $\# \H(P,R) \leq \# \H(P,S)$ for every finite poset $P$. Besides the structural investigation, there is also a need for an efficient toolbox for checking if 
$\# \H(P,R) \leq \# \H(P,S)$ holds for every finite poset $P$ or not, i.e., for checking if a strong Hom-scheme or even a strong G-scheme exists from $R$ to $S$. (The necessary conditions in \cite{aCampo_toappear_1} and the calculation rules and cancellation rules in \cite{aCampo_toappear_2} contribute to this toolbox.) The recent paper works on both aspects. Theorem \ref{theo_GschemeOnStrict} and Corollary \ref{coro_Lovaszrefinement} belong to the structural work, whereas the application-oriented Theorems \ref{theo_Transport} and \ref{theo_PQ_AT} are part of the toolbox.

\section{Preparation} \label{sec_preparation}

\subsection{Basics and Notation} \label{subsec_notation}

Let $X$ be a set. A reflexive, antisymmetric, and transitive relation $\leq \; \subseteq X \mytimes X$ is called a {\em partial order relation}, the pair $P = (X,\leq)$ is called a {\em partially ordered set} or simply a {\em poset}, and $X$ is called the {\em carrier} of $P$. As usual, we write $x \leq y$ for $(x,y) \in \, \leq$. For a poset $P = (X,\leq)$, ``$x < y$'' means ``$x \leq y$ and $x \not= y$''. For a subset $A \subseteq X$, the {\em poset induced on $A$} is defined as $P \vert_A \equiv ( A, \leq \cap \; ( A \mytimes A ) )$.

For a set $X$, the {\em diagonal (relation)} is defined as $\mysetdef{\Delta_X}{(x,x)}{x\in X}$, and $(X, \Delta_X)$ is called an {\em antichain}. A {\em chain} is characterized by $x \leq y$ or $y \leq x$ for all $x, y \in X$; up to isomorphism, there is only one chain for every underlying set. For a finite set $X$ of cardinality $k \in \myN_0$, we write $A_k$ for the antichain on $X$, and $C_k$ for the chain on $X$ (defined up to isomorphism).

\begin{figure} 
\begin{center}
\includegraphics[trim = 70 730 200 70, clip]{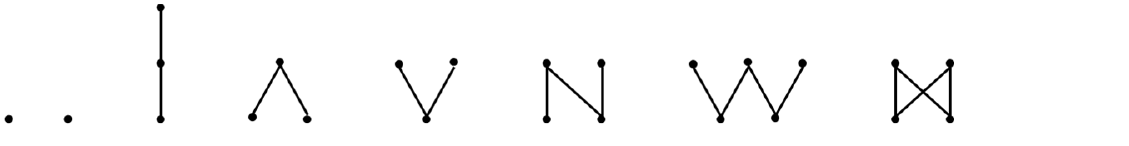}
\caption{\label{fig_basicPosets} The posets $A_2, C_3, \Lambda_3, V_3, N, W$ and $N^{(2)}$.}
\end{center}
\end{figure}

For posets $P_1 = (X_1, \leq_1)$ and $P_2 = (X_2, \leq_2)$ with $X_1 \cap X_2 = \emptyset$, their {\em direct sum} is the poset $P_1 + P_2 \equiv (X_1 \cup X_2, \leq_1 \cup \leq_2)$, and their {\em ordinal sum} is the poset $P_1 \oplus P_2 \equiv (X_1 \cup X_2, \leq_1 \cup \leq_2 \cup ( X_1 \mytimes X_2 ))$. For $k \in \myN$, we define $\Lambda_k \equiv A_{k-1} \oplus A_1$ (the bug with $k-1$ legs) and $V_k \equiv A_1 \oplus A_{k-1}$ (the $\Lambda_k$-bug turned onto its back). Furthermore, $N$ and $W$ are the posets with N- and W-shaped diagram, respectively, and $N^{(2)} \equiv A_2 \oplus A_2$ is the dubble-N. The diagrams of $A_2, C_3, \Lambda_3, V_3, N, W$, and $N^{(2)}$ are shown in Figure \ref{fig_basicPosets}.

Given posets $P = (X, \leq_P ) $ and $Q = (Y, \leq_Q)$, their product $P \mytimes Q \equiv (X \mytimes Y,$ $ \leq_{P \mytimes Q}) $ on $X \mytimes Y$ is defined by $(x_1,y_1) \leq_{P \mytimes Q} (x_2, y_2) $ iff $ x_1 \leq_P x_2$ and $y_1 \leq_Q y_2$. We write $P^k$ for $P \mytimes \ldots \mytimes P$ with $k$ factors $P$. A {\em binary word} is an element of $(C_2)^k$ with $0 < 1$ as carrier of $C_2$.

For posets $P = ( X, \leq_P ), Q = ( Y, \leq_Q )$, a mapping $\xi : X \rightarrow Y$ is called a {\em homomorphism}, iff $\xi(x) \leq_Q \xi(y)$ holds for all $x, y \in X$ with $x \leq_P y$. A homomorphism $\xi$ is called {\em strict} iff it fulfills additionally $x <_P y \Rightarrow \xi(x) <_Q \xi(y)$ for all $x, y \in X$, and it is called an {\em embedding} iff $\xi(x) \leq_Q \xi(y) \Rightarrow x \leq_P y$ for all $x, y \in X$. Finally, an embedding is called an {\em isomorphism} iff it is onto. $P \simeq Q$ indicates isomorphism. For posets $P$ and $Q$, we use the following symbols for homomorphism sets:
\begin{align*}
\H(P,Q) & \equiv \mysetdescr{ \xi : X \rightarrow Y }{ \xi \mytext{ is a homomorphism from} P \mytext{to} Q }, \\
\S(P,Q) & \equiv \mysetdescr{ \xi \in \H(P,Q) }{ \xi \mytext{is strict} }, \\
\S^o(P,Q) & \equiv \mysetdescr{ \xi \in \H(P,Q) }{ \xi \mytext{is strict and onto} }, \\
\Emb(P,Q) & \equiv \mysetdescr{ \xi \in \H(P,Q) }{ \xi \mytext{is an embedding} }, \\
\Aut(P) & \equiv \mysetdescr{ \xi \in \H(P,P) }{ \xi \mytext{is an isomorphism} }.
\end{align*}
We equip $\H(P,Q)$ with the ordinary pointwise partial order: for $\xi, \zeta \in \H(P,Q)$, we set $\xi \leq_{\H(P,Q)} \zeta$ iff $\xi(x) \leq_Q \zeta(x)$ for all $x \in P$.

$\mfP$ is the class of all finite posets, and the set $\mfP_r$ is a representation system of the non-isomorphic posets in $\mfP$.

We apply the terminology of homomorphisms also to relations which are not partial orders: for given sets $X$ and $Y$ and binary relations $R \subseteq X \times X$, $S \subseteq Y \times Y$, we call a mapping $f : X \rightarrow Y$ a homomorphism iff $(x,y) \in R$ implies $(f(x),f(y)) \in S$ for all $x, y \in X$, and we call a homomorphism strict, iff additionally $f(x) \not= f(y)$ holds for all $(x,y) \in R$ with $x \not= y$.

Let $X$ be a set and $R \subseteq X \mytimes X$ a binary relation on $X$. With $\R$ denoting the set of all transitive relations $S \subseteq X \mytimes X$ with $R \subseteq S$, the {\em transitive hull} $ T \equiv \cap \R$ of $R$ is the (set-theoretically) smallest transitive relation on $X$ containing $R$. If $R$ is reflexive, then also its transitive hull is reflexive; however, antisymmetry of $R$ is in general not preserved. $x T y$ is equivalent to the existence of $z_0, \ldots , z_L \in X$, $L \in \myN$, with $x = z_0, y = z_L$ and $z_{\ell-1} R z_\ell$ for all $1 \leq \ell \leq L$.

In order to avoid repetitions, we agree on that $X$ is always the carrier of the poset $P$, and that $Y$ is always the carrier of the poset $Q$. For a poset $P$, we use the notation $x \in P$ instead of $x \in X$, and for posets $P$ and $Q$, we write $\xi : P \rightarrow Q$ instead of $\xi : X \rightarrow Y$ for a homomorphism $\xi \in \H(P,Q)$.

In Section \ref{subsec_proving_RGS}, we need the following lemma:

\begin{lemma} \label{lemma_bijHom}
Let $P, Q \in \mfP$ and let $\xi \in \H(P,Q)$ be a bijective homomorphism. Then $\# ( \leq_P ) \; \leq \# ( \leq_Q )$ with equality iff $\xi$ is an isomorphism.
\end{lemma}
\BP We define the mapping $ \phi : \, \leq_P \rightarrow \leq_Q$ by $ (x,y)  \mapsto ( \xi(x), \xi(y) )$. It is one-to-one, and $\# ( \leq_R ) \, \leq \# ( \leq_S )$ is proven. Furthermore:
\begin{align*}
& \exists \; (x,y) \in \; \leq_Q \setminus \; \phi( \leq_P ) \\
\stackrel{\xi \mytext{bijective}}{\Leftrightarrow} & \exists \; (a,b) \in \; X^2 \; : \; (\xi(a), \xi(b) ) \in \; \leq_Q \; \mytext{and} \; ( a, b ) \, \notin \; \leq_P \\
\Leftrightarrow \quad & \xi \mytext{is not an embedding.}
\end{align*}

\EP

From the rich concept of downsets and upsets, we need simple notation only. Given a poset $P = ( X, \leq )$, we define for $x \in P$
\begin{align*}
\mydownarrow x & \equiv \mysetdescr{ y \in P }{ y \leq x }, \\
{\myodownarrow} x & \equiv \mysetdescr{ y \in P }{ y < x }, \\
\myuparrow x & \equiv \mysetdescr{ y \in P }{ x \leq y }, \\
{\myouparrow} x & \equiv \mysetdescr{ y \in P }{ x < y }.
\end{align*}
If required, we label the arrows with the relation they are referring to. A subset $A \subseteq X$ is called {\em convex} iff $ \myuparrow x \cap \mydownarrow y \subseteq A$ for every $x, y \in A$.

Additionally, we use the following notation from set theory:
\begin{align*}
\myNk{0} & \equiv  \emptyset, \\
\myNk{n} & \equiv  \{ 1, \ldots, n \} \mytext{for every} n \in \myN.
\end{align*}

For sets $X$ and $Y$, $\A(X,Y)$ is the set of mappings from $X$ to $Y$. $id_X \in \A(X,X)$ is the {\em identity mapping} defined by $x \mapsto x$. Let $f \in \A(X,Y)$ be a mapping. Our symbols for the {\em pre-image} of $B \subseteq Y$ and of $y \in Y$ are
\begin{align*}
\myurbild{f}(B) & \equiv \mysetdescr{ x \in X }{ f(x) \in B}, \\
\myurbild{f}(y) & \equiv \myurbild{f} (\{ y\} ).
\end{align*}
Furthermore, for any set $Z$ with $f(X) \subseteq Z$, we write $f \vert^Z$ for $ in_{f(X),Z} \circ f$, where $in_{f(X),Z}$ is the canonical inclusion from $f(X)$ into $Z$. For a one-to-one mapping $f \in \A(X,Y)$, the mapping $f^{-1} \in \A(f(X),X)$ is the inverse of $f$ on $f(X)$.

Finally, we use the {\em Cartesian product}. Let $\I$ be a non-empty set, and let $N_i$ be a non-empty set for every $i \in \I$. Then the Cartesian product of the sets $N_i, i \in \I$, is defined as
\begin{eqnarray*}
\prod_{i \in \I} N_i & \; \equiv \; & 
\mysetdescr{ f \in \A \big( \I, \bigcup_{i \in \I} N_i \big)}{ f(i) \in N_i \mytext{for all} i \in \I }.
\end{eqnarray*}

\subsection{Connectivity} \label{subsec_connectivity}

\begin{definition} \label{def_connected}
Let $P \in \mfP$, $A \subseteq P$, and $x, y \in A $. We say that $x$ and $y$ are {\em connected in $A$}, iff there are $z_0, z_1, \ldots , z_L \in A$, $L \in \myN_0$, with $x = z_0$, $y = z_L$ and $ z_{\ell-1} < z_\ell$ or $z_{\ell-1} > z_\ell$ for all $\ell \in \myNk{L}$. We call $z_0, \ldots , z_L$ a {\em zigzag line connecting $x$ and $y$}. We define for all $A \subseteq P$, $x \in A $
\begin{align} \nonumber
\gamma_A(x) & \equiv \mysetdescr{ y \in A }{ x \mytext{and} y \mytext{are connected in} A}.
\end{align}
\end{definition}
The relation ``connected in $A$'' is an equivalence relation on $A$ with partition $\mysetdescr{\gamma_A(a)}{a \in A}$.

The sets $\gamma_P(x), x \in P$, are called the {\em connectivity components} of $P$. Every poset is the direct sum of its connectivity components. A poset $P \in \mfP$ is {\em connected} iff $\gamma_P(x) = P$ for an $x \in P$ (the choice of $x \in P$ is arbitrary). We define
\begin{align*}
\mfP^c & \equiv \mysetdescr{ P \in \mfP }{ P \mytext{ connected } }, \\
\mfP_r^c & \equiv \mfP^c \cap \mfP_r.
\end{align*}
For a poset $P$ with connectivity components $Q_1, \ldots , Q_L$, $L \in \myN$, we have for every $R \in \mfP$
\begin{align*}
\H(P,R) & \; \simeq \; \H(Q_1, R) \times \cdots \times H(Q_L, R),
\end{align*}
and similar for $\S(P, R)$.

A subset $A \subseteq P$ is called {\em connected (in $P$)} iff $\gamma_A(x) = A$ for an (arbitrary) $x \in A$. For posets $P$ and $Q$, $P$ connected, the image $\xi(P)$ of $P$ under a homomorphism $\xi \in \H(P,Q)$ is connected in $Q$; in particular, $\xi(P)$ is a subset of a single connectivity component of $Q$. In consequence, if $P$ is a connected poset, then for all disjoint posets $R$ and $S$
\begin{align*}
\H(P, R + S) & \; \simeq \; \H(P,R) + \H(P,S),
\end{align*}
and similar for $\S(P, R+S)$ and $\Emb(P, R+S)$.

The following definition has already been used in earlier papers \cite{aCampo_2018,aCampo_toappear_1,aCampo_toappear_2}:
\begin{definition}
Let $P = (X, \leq)$ be a finite poset, let $Y$ be a set, and let $\xi \in \A(X,Y)$ be a mapping. We define for all $x \in X$
\begin{align*}
\gxix & \equiv \gamma_{\myurbild{\xi}( \xi(x) )}( x ).
\end{align*}
\end{definition}
$\gxix$ is thus the set of all points, $x$ is connected with in $\myurbild{\xi}( \xi(x) )$. As proven in \cite[Corollary 3]{aCampo_toappear_1}, a homomorphism $\xi$ is strict iff $\gxix = \{x\}$ for all $x \in P$.

\begin{lemma}[{\cite[Lemma 1]{aCampo_toappear_1}}] \label{lemma_Gxix_subsetnoteq_Gzetax}
Let $P, R, S \in \mfP$, let $\xi \in \H(P,R)$ and $\zeta \in \H(P,S)$ be homomorphisms, and let $\gxix \subseteq \gzex$ for an $x \in P$. Then $G_\xi(x) \subset G_\zeta(x)$ iff there are $a, b \in G_\zeta(x)$ with $a < b$ and $\xi(a) < \xi(b)$.
\end{lemma}

\section{Strong G-schemes} \label{sec_G_Schemes}

\subsection{Recapitulation} \label{subsec_Recapitulation}

\begin{figure} 
\begin{center}
\includegraphics[trim = 70 710 200 70, clip]{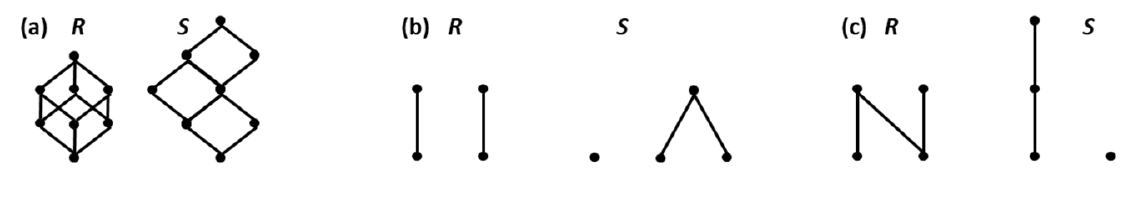}
\caption{\label{fig_Examples} Three pairs of posets $R$ and $S$ with $R \sqsubseteq_G S$.}
\end{center}
\end{figure}

In this section, we recall definitions and results contained in the previous papers \cite{aCampo_toappear_1,aCampo_toappear_2} about Hom-schemes.

\begin{definition} \label{def_RS_scheme}
Let $R, S \in \mfP$. We call a mapping
\begin{align*}
\rho & \in \prod_{P \in \mfP_r} \A( \H(P,R), \H(P,S) ) \quad \mytext{(Cartesian product)}
\end{align*}
a {\em Hom-scheme from $R$ to $S$}. We call a Hom-scheme $\rho$ from $R$ to $S$
\begin{itemize}
\item {\em strong} iff the mapping $\rho_P : \H(P,R) \rightarrow \H(P,S)$ is one-to-one for every $P \in \mfP_r$;
\item a {\em G-scheme} iff for every $P \in \mfP_r, \xi \in \H(P,R), x \in P$
\begin{align*}
G_{\rho_P(\xi)}(x) & = \gxix.
\end{align*}
\end{itemize}
If $P$ is fixed, we write $\rho(\xi)$ instead of $\rho_P(\xi)$.\end{definition}
We write 
\begin{equation*}
R \sqsubseteq S \; / \; R \sqsubseteq_G S
\end{equation*}
iff a strong Hom-scheme / a strong G-scheme exists from $R$ to $S$. The relations $\sqsubseteq$ and $\sqsubseteq_G$ define partial orders on $\mfP_r$ \cite[Theorems 2 and 3]{aCampo_toappear_1}.

If a one-to-one homomorphism $r : R \rightarrow S$ exists, we get a  (trivial) strong G-scheme from $R$ to $S$ by defining $\rho_P(\xi) \equiv r \circ \xi$ for all $P \in \mfP$, $ \xi \in \H(P,R)$. Non-trivial examples for posets $R$ and $S$ with $R \sqsubseteq_G S$ are shown in the Figures \ref{fig_Examples}, \ref{fig_NA1C3_WA1N2}, and \ref{fig_ConstrBspl}; more can be found in the Appendix of \cite{aCampo_toappear_1}. Postulating an additional regularity condition for the way, how a Hom-scheme maps the elements of $\H(P,R)$ to $\H(P,S)$, yields the concept of the {\em I-scheme}. An I-scheme is always a G-scheme; the structure theory of strong I-schemes is presented in \cite{aCampo_toappear_1}.

In the theory of strong G-schemes and strong I-schemes, the {\em EV-systems} play an important role:

\begin{definition} \label{def_EVsys}
Let $P$ be a poset. The {\em EV-system} $\E(P)$ of $P$ is defined as
\begin{align*}
\E(P) & \equiv \mysetdescr{ ( x, D, U ) }{ x \in P, D \subseteq  {\myodownarrow} x, U \subseteq {\myouparrow} x }.
\end{align*}
For $\fa \in \E(P)$, we refer to the three components of $\fa$ by $\fa_1, \fa_2$, and $\fa_3$. We equip $\E(P)$ with a relation: For all $\fa, \fb \in \E(P)$ we define
\begin{equation*}
\fa <_+ \fb \quad \equiv \quad \fa_1 \in \fb_2 \; \mytext{and} \; \fb_1 \in \fa_3,
\end{equation*}
and $\leq_+ \; \equiv \; <_+ \cup \; \Delta_{\E(P)}$. Additionally, we define for every $x \in P$
\begin{align*}
\E(P; x) & = \mysetdescr{ \fa \in \E(P) }{ \fa_1 = x }.
\end{align*}

\end{definition}

\begin{figure} 
\begin{center}
\includegraphics[trim = 80 630 185 70, clip]{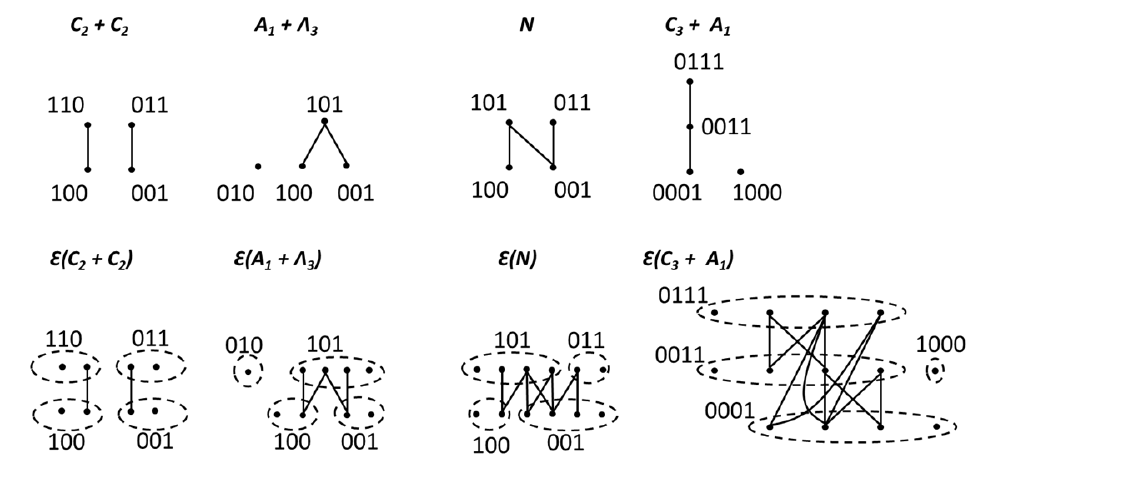}
\caption{\label{figure_EVSystems} The posets $C_2 + C_2, A_1 + \Lambda_3, N$, and $C_3 + A_1$ and their EV-systems. The points of the posets are labeled by binary words. In the EV-systems, the sets $\E(P;x)$ with $x \in P$ are encircled and labeled with the respective $x$.}
\end{center}
\end{figure}

Figure \ref{figure_EVSystems} shows the EV-systems of the posets $C_2 + C_2, A_1 + \Lambda_3, N$, and $C_3 + A_1$, and the EV-systems of $C_3$ and $N^{(2)}$ are shown in the Figures \ref{fig_N_C3} and \ref{fig_W_N2}. In all these figures, the sets $\E(P;x)$, $x \in P$, are encircled and labeled with $x$. For each point $\fa$ in the diagrams, $\fa_1$ is thus given by this label, and we get $\fa_2$ and $\fa_3$ by looking at the labels of the end points of lines starting in $\fa$ and going downwards and upwards, respectively.

A strict homomorphism $\eps : \E(R) \rightarrow \E(S)$ induces a G-scheme from $R$ to $S$ \cite[Proposition 2]{aCampo_toappear_1}, and if $\eps$ fulfills some additional conditions, the G-scheme is strong \cite[Proposition 3]{aCampo_toappear_1}. The relation $\leq_+$ is reflexive and antisymmetric, but in general not transitive. The mapping $\E(P) \rightarrow P$ defined by $\fa \mapsto \fa_1$ is a strict homomorphism, and the mapping $P \rightarrow \E(P)$ with $x \mapsto \left( x, {\myodownarrow} x, {\myouparrow} x \right)$ is an embedding. From the calculation rules for EV-systems presented in \cite{aCampo_toappear_2}, we need $\E(P+Q) = \E(P) + \E(Q)$ only.

In order to develop the structure theory of strong Hom-schemes, we have defined for every $P, Q \in \mfP$, $\xi \in \H(P,Q)$ a homomorphism $\alpha_\xi : P \rightarrow \E(Q)$ in \cite[Definition 5]{aCampo_toappear_1}. In the recent article, we will use these homomorphisms for {\em strict} homomorphisms $\xi$ only; we can thus use the following restricted and simplified definition:

\begin{definition} \label{def_axix}
Let $P, Q \in \mfP$, and let $\xi \in \H(P,Q)$ be a {\em strict} homomorphism. We define for every $x \in P$
\begin{align*}
\myapxix{P}{\xi}{x} & = \left( \xi(x), \xi( \myodownarrow x ), \xi( \myouparrow x ) \right).
\end{align*}
If $P$ is fixed, we write $\axix$ instead of $\myapxix{P}{\xi}{x}$.
\end{definition}
According to \cite[Lemma 3]{aCampo_toappear_1}, $\xi \in \S(P,Q)$ implies $\alpha_\xi \in \S(P, \E(Q))$.

\subsection{Strong G-schemes and strict homomorphisms} \label{subsec_GschemeOnStrict}

In this section, we prove

\begin{theorem} \label{theo_GschemeOnStrict}
Let $R, S \in \mfP$. Equivalent are
\begin{align}
R & \sqsubseteq_G S; \label{r_homG_S} \\
\# \S(P,R) & \leq \# \S(P,S) \; \; \mytext{for all} \; P \in \mfP;  \label{SPR_leq_SPS} \\
\# \S(Q,R) & \leq \# \S(Q,S) \; \; \mytext{for all} \; Q \in \mfP^c. \label{SPR_leq_SPS_connected}
\end{align}
\end{theorem}
The proof of $\eqref{r_homG_S} \Rightarrow \eqref{SPR_leq_SPS} \Leftrightarrow \eqref{SPR_leq_SPS_connected}$ is simple. Assume that \eqref{r_homG_S} holds. As stated in Section \ref{subsec_connectivity}, a homomorphism $\xi \in \H(P,Q)$ is strict iff $\gxix = \{ x \}$ for every $x \in P$. For a G-scheme $\rho$ from $R$ to $S$, we have $\grxix = \gxix$ for every $\arie$, thus $\rho_P( \S(P,R) ) \subseteq S(P,S)$ for every $P \in \mfP_r$, and \eqref{SPR_leq_SPS} follows if $\rho$ is strong.

$\eqref{SPR_leq_SPS} \Rightarrow \eqref{SPR_leq_SPS_connected}$ is trivial. Assume that \eqref{SPR_leq_SPS_connected} holds, and let $P \in \mfP$ with connectivity components $Q_1, \ldots , Q_L$. Then
\begin{align*}
\# \S(P,R)  \; & = \; \prod_{\ell=1}^L \# \S(Q_\ell, R )  \; \leq \;
\prod_{\ell=1}^L \# \S(Q_\ell, S ) \; = \; \# \S(P,S).
\end{align*}

It remains to prove $\eqref{SPR_leq_SPS} \Rightarrow \eqref{r_homG_S}$. For this proof, we modify the well-known mathematical approach to let a mapping $f : A \rightarrow B$ factorize over the set $U(f) \equiv \mysetdescr{ \myurbild{f}(f(a)) }{ a \in A}$ of its pre-images: $f = \iota_f \circ \pi_f$ with canonical mappings $\pi_f : A \rightarrow U(f)$, $\iota_f : U(f) \rightarrow B$. In this section, we let a homomorphism factorize over a {\em refinement} of the set of its pre-images, a refinement consisting of the {\em connected} sets $\gxix$:

\begin{definition}
Let $P, Q \in \mfP$. We define for every $ \xi \in \H(P,Q)$ the set
\begin{align*}
\Gxi \; & \equiv \; \mysetdescr{ \gxix }{ x \in P },
\end{align*}
and the mappings
\begin{align*}
\pi_\xi : P & \rightarrow \Gxi \\
x & \mapsto \gxix, \\
\iota_\xi : \Gxi & \rightarrow Q \\
\gxix & \mapsto \xi(x).
\end{align*}
Additionally, we define a relation $\preceq_0$ on $\Gxi$ by setting for all $\fa, \fb \in \Gxi$
\begin{align*}
\fa \preceq_0 \fb \; & \equiv \; \exists \; a \in \fa, b \in \fb \; : \; a \leq b.
\end{align*}
$\preceq$ is the transitive hull of $\preceq_0$. \end{definition}

$\Gxi$ is a partition of $X$ consisting of connected sets, and $\iota_\xi$ is a well-defined mapping because $\xi$ is constant on every $\fc \in \Gxi$: $\xi(x) = \xi(y)$ for all $x, y \in \fc$. Obviously, $\xi = \iota_\xi \circ \pi_\xi$ for all $\xi \in \H(P,Q)$. Furthermore:
\begin{lemma} \label{lemma_iota}
Let $P, Q \in \mfP$ and $ \xi \in \H(P,Q)$, and let $\fc_0, \ldots , \fc_I \in \Gxi$, $I \in \myN$, with $\fc_{i-1} \preceq_0 \fc_i$ for every $i \in \myNk{I}$. Then $\iota_\xi( \fc_0 ) \leq \iota_\xi( \fc_I )$, and $\iota_\xi( \fc_0 ) = \iota_\xi( \fc_I )$ implies $\fc_0 = \fc_1 = \ldots = \fc_I$.
\end{lemma}
\BP $\fc_{i-1} \preceq_0 \fc_i$ for every $i \in \myNk{I}$ is equivalent to the existence of $x_0^+ \in \fc_0$, $x_i^-, x_i^+ \in \fc_i$ for every $i \in \myNk{{I-1}}$, and $x_I^- \in \fc_I$ with $x_{i-1}^+ \leq x_i^-$ for every $i \in \myNk{I}$. Because the sets $\fc_i$ are all connected in $P$, we conclude that also the set $V \equiv \bigcup_{i=0}^I \fc_i$ is connected in $P$.

We have $\xi( x_{i-1}^+ ) \leq \xi( x_i^- )$ for every $i \in \myNk{I}$, and because $\xi$ is constant on every set $\fc_i$, we have $\xi( x_i^- ) = \xi( x_i^+ )$ for every $i \in \myNk{{I-1}}$. Together this yields the following chain of inequalities and equations:
\begin{align*}
\xi(x_0^+ ) \; & \leq \; \xi( x_1^- ) \; = \; \xi(x_1^+) \; \leq \; \xi( x_2^- ) \; = \; \xi(x_2^+) \\
& \leq \; \ldots \; \leq \; \xi( x_{I-1}^- ) \; = \; \xi(x_{I-1}^+) \; \leq \; \xi( x_I^- ).
\end{align*}
We thus have $\iota_\xi( \fc_0 ) = \xi(x_0^+ ) \leq \xi( x_I^- ) = \iota_\xi( \fc_I )$. In the case of $\iota_\xi( \fc_0 ) = \iota_\xi( \fc_I )$, we get $ \xi(x_0^+ )  = \xi( x_I^- )$, and we conclude that all inequalities in the chain above are in fact equations. $\xi$ is thus constant on the connected set $V$, hence $\fc_i = \mygxix{\xi}{x_i^-} = \gamma_{\myurbild{\xi}(\xi(x_i^-))}(x_i^-) \supseteq V$ for every $i \in \myNk{{I}}$, and similarly $\fc_0 = \mygxix{\xi}{x_0^+} \supseteq V$, which yields $\fc_i = V$ for every $i \in \myNk{{I}} \cup \{ 0 \}$.

\EP

\begin{corollary} \label{coro_preceq_pi_iota}
Let $P, Q \in \mfP$ and $ \xi \in \H(P,Q)$. Then $\preceq$ is a partial order on $\Gxi$, $\pi_\xi : P \rightarrow \Gxi$ is a homomorphism, and $\iota_\xi : \Gxi \rightarrow Q$ is a strict homomorphism.
\end{corollary}
\BP $\preceq$ is a reflexive and transitive relation on $\Gxi$. We have to show that $\preceq$ is antisymmetric.

Let $\fa, \fb \in \Gxi$ with $\fa \preceq \fb$ and $\fb \preceq \fa$. There exist $\fc_0, \ldots , \fc_I, \fc_{I+1}, \ldots ,$ $ \fc_{I+J} \in \Gxi$, $I, J \in \myN$, with $\fa = \fc_0 = \fc_{I+J}$, $\fb = \fc_I$, and $\fc_{i-1} \preceq_0 \fc_i$ for every $i \in \myNk{{I+J}}$. Because of $\fc_0 = \fa = \fc_{I+J}$ we have $\iota_\xi( \fc_0 ) = \iota_\xi( \fc_{I+J} )$. Lemma \ref{lemma_iota} delivers $\fc_0 = \fc_1 = \ldots = \fc_{I+J}$, and $\fa = \fb$ is shown.

Now it is trivial that $\pi_\xi$ is a homomorphism. The first statement about $\iota_\xi$ in Lemma \ref{lemma_iota} tells us that $\iota_\xi$ is a homomorphism, and the second one yields that it is a strict one.

\EP

\begin{lemma} \label{lemma_Gscheme_Gamma}
For every $P, Q \in \mfP$, $ \xi \in \H(P,Q)$, we define for every $T \in \mfP$
\begin{align*}
\Gamma_{P,T}(\xi) & \equiv \mysetdescr{ \zeta \in \H(P,T) }{ \Gzet = \Gxi }.
\end{align*}
Then, for $R, S \in \mfP$, the relation $R \sqsubseteq_G S$ is equivalent to
\begin{align*}
\# \Gamma_{P,R}(\xi) \leq \# \Gamma_{P,S}(\xi) \; \; \mytext{for all} \; P \in \mfP, \xi \in \H(P,R).
\end{align*}
\end{lemma}
\BP Let $R, S \in \mfP$, $P \in \mfP$, $\xi \in \H(P,R)$, $\zeta \in \H(P,S)$. Because $\Gxi$ and $\Gzet$ are both partitions of $X$, there is for every $x \in X$ a unique $\fa \in \Gxi$ and a unique $\fb \in \Gzet$ with $x \in \fa = \gxix$ and $x \in \fb = \gzex$. $\Gzet = \Gxi$ is thus equivalent to $\gxix = \gzex$ for all $x \in P$. A G-scheme $\rho$ from $R$ to $S$ maps thus $\Gamma_{P,R}(\xi)$ to $\Gamma_{P,S}(\xi)$ for every $P \in \mfP_r$, $\xi \in \H(P,R)$, and if it is strong, it does so one-to-one. On the other hand, $\# \Gamma_{P,R}(\xi) \leq \# \Gamma_{P,S}(\xi)$ for every $P \in \mfP, \xi \in \H(P,R)$, gives trivially raise to a strong G-scheme because for every $P \in \mfP_r$, the set $\mysetdescr{ \Gamma_{P,R}(\xi) }{ \xi \in \H(P,R) }$ is a partition of $\H(P,R)$, and the set $\mysetdescr{ \Gamma_{P,S}(\xi) }{ \xi \in \H(P,R) }$ is a collection of disjoint subsets of $\H(P,S)$.

\EP

\begin{corollary} \label{coro_Gamma_pi_iota}
Let $P, Q \in \mfP$ and $\xi \in \H(P,Q)$. Then for every $T \in \mfP$
\begin{align} \label{eq_pi_zeta}
\pi_\zeta & = \pi_\xi \quad \mytext{for all} \; \zeta \in \Gamma_{P,T}(\xi),
\end{align}
and
\begin{align} \label{equiv_pi_zeta}
\zeta_1 \not= \zeta_2 \; & \Leftrightarrow \; \iota_{\zeta_1} \not= \iota_{\zeta_2} \quad \mytext{for all} \; \zeta_1, \zeta_2 \in \Gamma_{P,T}(\xi).
\end{align}
\end{corollary}
\BP Let $\zeta \in \Gamma_{P,T}(\xi)$. According to Corollary \ref{coro_preceq_pi_iota}, $\pi_\xi$ and $\pi_\zeta$ are  both elements of $\H(P, \Gxi)$. Let $x \in P$. Because $\Gxi = \Gzet$ is a partition of $X$, there is a unique $\fa \in \Gxi = \Gzet$ with $x \in \fa$. We conclude $\pi_\xi(x) = \gxix = \fa = \gzex = \pi_\zeta(x)$, and \eqref{eq_pi_zeta} is shown, because $x \in P$ was arbitrary.

For $\zeta_1, \zeta_2 \in \Gamma_{P,T}(\xi)$, Corollary \ref{coro_preceq_pi_iota} yields $\iota_{\zeta_1}, \iota_{\zeta_2} \in \S( \Gxi, T )$. The conclusion $\zeta_1 = \zeta_2 \Rightarrow \iota_{\zeta_1} = \iota_{\zeta_2}$ is trivial. Let $\zeta_1 \not= \zeta_2$. Then $\iota_{\zeta_1} \circ \pi_{\zeta_1} = \zeta_1 \not= \zeta_2 = \iota_{\zeta_2} \circ \pi_{\zeta_2}$, and \eqref{eq_pi_zeta} delivers $\iota_{\zeta_1} \not= \iota_{\zeta_2}$

\EP

\begin{lemma} \label{lemma_Gamma_SPGxiQ}
For $P, Q, T \in \mfP$, $\xi \in \H(P,Q)$, we have
$\# \Gamma_{P,T}(\xi) = \# \S( \Gxi, T )$.
\end{lemma}
\BP Let $\J \equiv \mysetdescr{ \iota_\zeta }{ \zeta \in \Gamma_{P,T}(\xi) }$. The last statement in Corollary \ref{coro_preceq_pi_iota} yields $\J \subseteq \S(\Gxi,T)$, and \eqref{equiv_pi_zeta} delivers $ \# \Gamma_{P,T}(\xi) = \# \J$. We conclude $\# \Gamma_{P,T}(\xi) \leq \# \S( \Gxi, T )$.

Let now $\sigma_1, \sigma_2 \in \S( \Gxi, T )$ with $\sigma_1 \not= \sigma_2$. With $\zeta_1 \equiv \sigma_1 \circ \pi_\xi$, $\zeta_2 \equiv \sigma_2 \circ \pi_\xi$ we have $\zeta_1, \zeta_2 \in \H(P,T)$ because of $\pi_\xi \in \H(P,\Gxi)$ (Corollary \ref{coro_preceq_pi_iota}).

Let $i \in \myNk{2}$ be fixed. We want to show $\zeta_i \in \Gamma_{P,T}(\xi)$. For $x \in P$, the set $\gxix$ contains $x$ and is connected in $P$, and the mapping $\zeta_i$ is constant on $\pi_\xi(x) = \gxix$. Therefore, $\gxix \subseteq \mygxix{\zeta_i}{x}$. In the case of ``$\subset$'', Lemma \ref{lemma_Gxix_subsetnoteq_Gzetax} delivers
$a, b \in \mygxix{\zeta_i}{x}$ with $a < b$ and $\xi(a) < \xi(b)$. But $\xi(a) < \xi(b)$ means $\mygxix{\xi}{a} \not= \mygxix{\xi}{b}$, and $a < b$ means $\mygxix{\xi}{a} = \pi_\xi(a) \preceq \pi_\xi(b) = \mygxix{\xi}{b}$, which together yields $\mygxix{\xi}{a} \prec \mygxix{\xi}{b}$. Because $\sigma_i$ is strict, we have
\begin{displaymath}
\zeta_i(a) = \sigma_i ( \pi_\xi( a ) ) = \sigma_i( \mygxix{\xi}{a} ) < \sigma_i( \mygxix{\xi}{b} ) = \sigma_i ( \pi_\xi( b ) ) = \zeta_i(b)
\end{displaymath}
in contradiction to $a, b \in \mygxix{\zeta_i}{x} \subseteq \myurbild{\zeta_i}(\zeta_i(x))$. Therefore, $\gxix = \mygxix{\zeta_i}{x}$, thus $\zeta_i \in \Gamma_{P,T}(\xi)$, because $x \in P$ was arbitrary.

We have $\iota_{\zeta_1} = \sigma_1$ and $\iota_{\zeta_2} = \sigma_2$. Equivalence \eqref{equiv_pi_zeta} yields $\zeta_1 \not= \zeta_2$, and $\# \S( \Gxi, T ) \leq \# \Gamma_{P,T}(\xi)$ is shown.

\EP

Now we can prove the implication $\eqref{SPR_leq_SPS} \Rightarrow \eqref{r_homG_S}$ in Theorem \ref{theo_GschemeOnStrict}. Let $R, S \in \mfP$ with $\# \S(P,R) \leq \# \S(P,S)$ for all $P \in \mfP$. According to Lemma \ref{lemma_Gscheme_Gamma}, we have to show $\# \Gamma_{P,R}(\xi) \leq \# \Gamma_{P,S}(\xi)$ for all $P \in \mfP, \xi \in \H(P,R)$. And that is easily done. Let $P \in \mfP$, $\xi \in \H(P,R)$. Applying Lemma \ref{lemma_Gamma_SPGxiQ} with $Q = T = R$ yields $\# \Gamma_{P,R}(\xi) = \# \S( \Gxi, R )$, and using Lemma \ref{lemma_Gamma_SPGxiQ} with $Q = R, T = S$ results in $\# \Gamma_{P,S}(\xi) = \# \S( \Gxi, S )$. Because \eqref{SPR_leq_SPS} delivers $ \# \S( \Gxi, R ) \leq \# \S(\Gxi,S)$, the inequality $\# \Gamma_{P,R}(\xi) \leq \# \Gamma_{P,S}(\xi)$ is shown.

In the Figures \ref{fig_NA1C3_WA1N2} and \ref{fig_ConstrBspl} in Section \ref{sec_application}, pairs of posets $R$ and $S$ are shown with $R \sqsubseteq_G S$. In all pairs, we have $S = A + T$ with an antichain $A$ and a connected poset $T$. The equivalence $\eqref{r_homG_S} \Leftrightarrow \eqref{SPR_leq_SPS_connected}$ in Theorem \ref{theo_GschemeOnStrict} sheds an interesting light on the roles of $A$ and $T$ in $R \sqsubseteq_G S$. For all $Q \in \mfP^c$ with $\# Q \geq 2$, $\# \S(Q,R) \leq \# \S(Q,A+T)$ is equivalent to $\# \S(Q,R) \leq \# \S(Q,T)$. It is thus $T$ which bears the main burden of $R \sqsubseteq_G S$, whereas the only role of $A$ is to ensure $\# \S(A_1,R) \leq \# \S(A_1,S)$ by adding sufficiently many points.

As a first application of Theorem \ref{theo_GschemeOnStrict}, we refine a result of Lov\'{a}sz \cite{Lovasz_1967}:

\begin{corollary} \label{coro_Lovaszrefinement}
Let $R, S$ be finite posets. Then
\begin{align*}
\# \S(P,R) = \# \S(P,S) \; \mytext{for every} \; P \in \mfP \quad \Leftrightarrow \quad R \simeq S.
\end{align*}
\end{corollary}
\BP $\# \S(P,R) = \# \S(P,S)$ for every $P \in \mfP$ means $R \sqsubseteq_G S$ and $S \sqsubseteq_G R$ according to Theorem \ref{theo_GschemeOnStrict}. Because $\sqsubseteq_G$ is a partial order on $\mfP_r$, we get $R \simeq S$.

\EP

Furthermore, we simplify the conditions in \cite[Proposition 3]{aCampo_toappear_1}:

\begin{proposition} \label{prop_eta_invers}
Let $\eps : \E(R) \rightarrow \E(S)$ be a strict homomorphism between the EV-systems of the posets $R, S \in \mfP$. For every poset $P \in \mfP^c_r$ and for every homomorphism $\xi \in \S(P,R)$, we define for every $x \in P$
\begin{align*}
\myepxix{P}{\xi}{x} & \equiv \eps( \myapxix{P}{\xi}{x} )_1,
\end{align*}
Assume that $\eps( \fa ) = \eps( \fb )$ implies $\fa_1 = \fb_1$ for every $\fa, \fb \in \E(R)$. Assume additionally, that for every $P \in \mfP^c_r$, $\xi \in \S(P,R)$, $x \in P$
\begin{align*}
\myapepxix{P}{\xi}{x} \in \eps( \E(R)) & \quad \Rightarrow \quad  \myapepxix{P}{\xi}{x} \in \eps( \E( R; \xi(x) ) ), \\
\myapepxix{P}{\xi}{x} \notin \eps( \E(R)) & \quad \Rightarrow \quad  \xi(x) \in Z \setminus Z^+,
\end{align*}
where $Z^+$ is a subset of the carrier $Z$ of $R$ with $ \# ( Z \setminus Z^+ ) \leq 1$. Then $R \sqsubseteq_G S$ with
\begin{equation} \label{eq_eta_invers}
\myurbild{\xi}(v) \quad = \quad \bigcup_{\fa \in \E(R;v)}\myurbild{ \apepxi }( \eps( \fa ) )
\end{equation}
for every $v \in Z^+$, $P \in \mfP^c_r$, $\xi \in \S(P,R)$.
\end{proposition}
\BP As mentioned after Definition \ref{def_EVsys}, the mapping $\E(S) \rightarrow S$ with $ \fa \mapsto \fa_1$ is a strict homomorphism. Because of $\alpha_\xi \in \S(P,\E(R))$ for every $P \in \mfP^c_r$ and every $\xi \in \S(P,R)$, $\eta_P$ is a well-defined mapping from $\S(P,R)$ to $\S(P,S)$ for every $P \in \mfP^c_r$. Now equation \eqref{eq_eta_invers} is proven exactly as in the proof of \cite[Proposition 3]{aCampo_toappear_1}, and just as there we conclude that $\eta_P$ maps $\S(P,R)$ one-to-one to $\S(P,S)$ for every $P \in \mfP^c_r$. Now apply Theorem \ref{theo_GschemeOnStrict}.

\EP

\section{A sufficient condition for $R \sqsubseteq_G S$} \label{sec_suffCond}

We have seen in Theorem \ref{theo_GschemeOnStrict} that it is the relation between the cardinalities of the sets $\S(P,R)$ and $\S(P,S)$ over all $P \in \mfP$ which is crucial for the existence of a  strong G-scheme from $R$ to $S$. If we want to apply this result (or Proposition \ref{prop_eta_invers}) in proving $R \sqsubseteq_G S$, we are still dealing with an infinite number of posets $P$. Building on Theorem \ref{theo_GschemeOnStrict}, we develop in this section a sufficient condition for $R \sqsubseteq_G S$ which needs the inspection of a {\em finite} set of posets only. As starting point, we recall a factorization and results due to Lov\'{a}sz \cite{Lovasz_1967}.

\begin{definition} \label{def_IQPT}
For all $Q, P, T \in \mfP$, we define
\begin{align*}
\mfE(T) & \equiv \mysetdescr{ E \in \mfP_r }{ E \; \mytext{can be embedded in} T }, \\
\mfE^c(T) & \equiv \mfE(T) \cap \mfP_r^c, \\
\I_Q(P,T) & \equiv \mysetdescr{ \xi \in \S(P,T) }{ \xi(P) \simeq Q }.
\end{align*}
\end{definition}
For $P, T \in \mfP$, the sets $\I_Q(P,T), Q \in \mfE(T)$, form a partition of $\S(P,T)$ \cite{Lovasz_1967}:
\begin{align} \label{eq_summe_bQPT}
\# \S(P,T) & = \sum_{ \in \mfE(T)} \# \I_Q(P,T).
\end{align}
If $P$ is connected, then $\I_Q(P,T) = \emptyset$ for every non-connected poset $Q \in \mfP$, hence
\begin{align} \label{eq_summe_bQPT_conn}
\# \S(P,T) & = \sum_{Q \in \mfE^c(T)} \# \I_Q(P,T).
\end{align}

Let $P, T \in \mfP$. For every $\xi \in \S(P,T)$, we have $\xi \in \I_{\xi(P)}(P,T)$. Moreover, we can factorize $\xi$ in the following way \cite{Lovasz_1967}. There exists a $Q \in \mfE(T)$ with $Q \simeq \xi(P)$, a surjective strict homomorphism $\sigma : P \rightarrow Q$, and an embedding $\epsilon : Q \rightarrow T$ with $\xi = \epsilon \circ \tau$. Now we define an equivalence relation on $\S^o(P,Q)$ by setting for all $\sigma, \sigma' \in \S^o(P,Q)$
\begin{displaymath}
\sigma \sim \sigma' \; \equiv \; \exists \, \beta \in \Aut(Q) \, : \sigma' = \beta \circ \sigma.
\end{displaymath}
With $\S_r^o(P,Q)$ denoting a representation system of the partition of $\S^o(P,Q)$ induced by this equivalence relation, we have for all $P, T \in \mfP$ for every $Q \in \mfE(T)$ \cite[Equation (7)]{Lovasz_1967}
\begin{align} \label{eq_formel_IQPT_1}
\# \I_Q(P,T) & = \# \S_r^o(P,Q) \cdot \# \Emb( Q, T ) \\
& = \# \S^o(P,Q) \cdot \frac{ \# \Emb( Q, T )}{ \# \Aut(Q) }. \label{eq_formel_IQPT_2}
\end{align}

Via Equation \eqref{eq_summe_bQPT_conn}, Equation \eqref{eq_formel_IQPT_1} is useful for the calculation of $\# \S(P,T)$ for not too complicated posets $P$ and $T$, $P$ connected. For those posets, it is easy to enumerate $\mfE^c(T)$ and to calculate $\# \S_r^o(P,Q) $ and $ \# \Emb( Q, T ) $ for the posets $Q \in \mfE^c(T)$. Two examples are shown in Figure \ref{fig_MatrixBspl} in matrix notation.

\setlength{\arraycolsep}{0.5mm}
\begin{figure}
$
\left( \begin{array}{ccccccc}
& A_1 & C_2 & V_3 & \Lambda_3 & N & C_3 \\
A_1 & 1 &&&&& \\
C_2 & & 1 &&&& \\
V_3 & & 1 & 1 &&& 2\\
\Lambda_3 & & 1 && 1 &&  2\\
N & & 1 & 1 & 1 & 1 & 5 \\
C_3 &&&&&& 1
\end{array} \right)
\left( \begin{array}{cc}
N & A_1 + C_3 \\
4 & 4 \\
3 & 3 \\
2 & \\
2 & \\
1 & \\
& 1 
\end{array} \right)
\quad = \quad
\left( \begin{array}{cc}
N & A_1 + C_3 \\
4 & 4 \\
3 & 3 \\
5 & 5 \\
5 & 5 \\
8 & 8 \\
& 1 
\end{array}
\right)
$
\\ \\ \\
$
\left( \begin{array}{cccccccc}
& A_1 & C_2 & V_3 & \Lambda_3 & N & W & N^{(2)} \\
A_1 & 1 &&&&&& \\
C_2 && 1 &&&&& \\
V_3 && 1 & 1 &&&& \\
\Lambda_3 && 1 && 1 &&& \\
N && 1 & 1 & 1 & 1 && 1 \\
W && 1 & 3 & 1 & 2 & 1 & 3 \\
N^{(2)} && 1 & 1 & 1 &&& 1
\end{array} \right)
\left( \begin{array}{cc}
W & A_1 + N^{(2)} \\
5 & 5 \\
4 & 4 \\
4 & 4 \\
2 & 4 \\
2 & \\
2 & \\
& 4 
\end{array} \right)
\quad = \quad
\left( \begin{array}{cc}
W & A_1 + N^{(2)} \\
5 & 5 \\
4 & 4 \\
8 & 8 \\
6 & 8 \\
12 & 16 \\
24 & 32 \\
10 & 16 
\end{array}
\right)
$
\caption{\label{fig_MatrixBspl} For the posets $R = N$, $S = A_1 + C_3$ in the upper part and for the posets $R = W$, $S = A_1 + N^{(2)}$ in the lower part, the following matrices are shown (from left to right, zeros are omitted):
$\left( \# \S_r^o(Q, Q') \right) $ with $Q, Q' \in \mfE^c(R) \cup \mfE^c(S)$, 
$\left( \# \Emb(Q, T) \right) $ with $Q \in \mfE^c(R) \cup \mfE^c(S), T \in \{ R, S \}$, and
$\left( \# S(Q, T) \right) $ with $Q \in \mfE^c(R) \cup \mfE^c(S), T \in \{ R, S \}$.
}
\end{figure}

The next corollary is a direct consequence of \eqref{eq_summe_bQPT_conn}, \eqref{eq_formel_IQPT_1}, and Theorem \ref{theo_GschemeOnStrict}:

\begin{corollary} \label{coro_applic_sum_of_soe}
Assume that $R, S \in \mfP$ are posets with $\mfE^c(R) \subseteq \mfE^c(S)$ and 
$ \# \Emb(Q,R) \leq \# \Emb(Q,S)$ for every $Q \in \mfE^c(R)$. Then $R \sqsubseteq_G S$.
\end{corollary}

Using this corollary, it is easily seen that $R \sqsubseteq_G S$ holds for the posets in Figure \ref{fig_Examples}(b). Furthermore, we will use it in Section \ref{subsec_constr_Gschemes} for the construction of strong G-schemes between posets of a special type. In the recent section, the corollary is the starting point of a generalization.

We can regard Corollary \ref{coro_applic_sum_of_soe} as a (trivial) example of a structured mapping of homomorphism sets. Let $Q \in \mfE^c(R)$ and assume that there exists a $Q' \in \mfE^c(S)$ with a one-to-one mapping of $\S^o(P, Q )$ into $\S^o(P, Q')$ ($Q' = Q$ in Corollary \ref{coro_applic_sum_of_soe}). In the case of $\frac{ \# \Emb(Q,R) }{ \# \Aut(Q) } \leq \frac{ \# \Emb(Q',S) }{ \# \Aut(Q') }$, Equation \eqref{eq_formel_IQPT_2} yields the existence of a one-to-one mapping from $\I_Q(P,R)$ to $\I_{Q'}(P,S)$.

The idea how to map homomorphism sets in a structured way is simple:

\begin{definition} \label{def_distributing}
Let $Q, Q' \in \mfP$, and let $\tau \in \S^o(Q, Q')$. We define for every $P \in \mfP^c$,
\begin{align*}
\T_\tau(P,Q) & \equiv \mysetdescr{ \tau \circ \sigma }{ \sigma \in \S^o(P,Q) } \; \subseteq \; \S^o(P,Q').
\end{align*}
We call $\tau$ {\em distributing}, iff the mapping $\S^o(P,Q) \rightarrow \T_\tau(P, Q)$, $\xi \mapsto \tau \circ \xi$, is one-to-one (i.e., bijective) for every $P \in \mfP^c$.	
\end{definition}

Let $P, Q, Q' \in \mfP$, $\sigma \in \S(P,Q)$, and $\tau \in \S(Q,Q')$. Then for all $x \in P$
\begin{align*}
\atsX & = \left( \tau(\sigma(x)), \tau(\sigma( \myodownarrow x ) ), \tau(\sigma( \myouparrow x ) ) \right), \\
\atSX & = \left( \tau(\sigma(x)), \tau( \myodownarrow \sigma( x ) ), \tau( \myouparrow \sigma( x ) ) \right).
\end{align*}
Because $\sigma$ is a strict homomorphism, we have $\sigma( \myodownarrow x ) \subseteq \, \myodownarrow \sigma(x)$ and $\sigma( \myouparrow x ) \subseteq \, \myouparrow \sigma(x)$, hence $\atsX_2 \subseteq \atSX_2$ and $\atsX_3 \subseteq \atSX_3$. Part (a) of the following lemma refers to this fact.

\begin{lemma} \label{lemma_trennbedingung}
Let $Q, Q' \in \mfP$ and $\tau \in \S^o(Q, Q')$ with $\alpha_\tau$ one-to-one.
\begin{itemize}
\item[(a)] Assume, that for every $P \in \mfP^c$, $\sigma \in \S^o(P,Q)$, $x \in P$, the point $\atSX$ is the only $\fa \in \alpha_\tau(Q)$ with $\atsX_1 = \fa_1$, $\atsX_2 \subseteq \fa_2$, and $\atsX_3 \subseteq \fa_3$. Then $\tau$ is distributing.
\item[(b)] The conditions in (a) are fulfilled (and $\tau$ is distributing) if for all $v, w \in Q$
\begin{align} \label{trennbedingung_2}
v \not= w \; \mytext{and} \; \tau(v) = \tau(w) \; \Rightarrow & \; \quad \, \quad \alpha_\tau(v)_2 \cap \alpha_\tau(w)_2 = \emptyset \\
& \; \mytext{and} \; \alpha_\tau(v)_3 \cap \alpha_\tau(w)_3 = \emptyset. \nonumber
\end{align}
\end{itemize}
\end{lemma}
\BP (a) Let $P \in \mfP^c$, $\sigma \in \S^o(P,Q)$, and $x \in P$. Because $\atSX$ is the only $\fa \in \alpha_\tau(Q)$ with $\atsX_1 = \fa_1$, $\atsX_2 \subseteq \fa_2$, and $\atsX_3 \subseteq \fa_3$, we can identify $\atSX$ by means of $\atsX$. In this way, we get a mapping
\begin{align} \label{def_phi_tau}
\phi_\tau : \mysetdescr{ \myapxix{P}{\theta}{x}}{ P \in \mfP^c, \theta \in \T_\tau(P,Q), x \in P } & \rightarrow \alpha_\tau(Q), \\
\atsX & \mapsto \atSX \nonumber.
\end{align}
Because $\alpha_\tau$ is one-to-one, we have $\alpha_\tau^{-1} \circ \phi_\tau \circ \alpha_{\tau \circ \sigma} = \sigma$ for all $\sigma \in \S^o(P,Q)$, and $\tau$ must be distributing.

(b) Let $P \in \mfP^c$, $\sigma \in \S^o(P,Q)$, and $x \in P$. We know already $\atsX_1 = \atSX_1$, $\atsX_2 \subseteq \atSX_2$, and $\atsX_3 \subseteq \atSX_3$

Let $\fa \in \alpha_\tau(Q)$ with $\atsX_1 = \fa_1$, $\atsX_2 \subseteq \fa_2$ and $\atsX_3 \subseteq \fa_3$. Because $\sigma$ is onto, there exists a $y \in P$ with $\fa = \atSY$. In particular, we have $\tau(\sigma(x)) = \atsX_1 = \fa_1 = \atSY_1 = \tau(\sigma(y))$. Assume $\sigma( y ) \not= \sigma(x)$. Then
\begin{align*}
\atsX_2 & \subseteq \atSX_2 \cap \atSY_2 \stackrel{\eqref{trennbedingung_2}}{=} \emptyset
\end{align*}
and similarly $\atsX_3 = \emptyset$. In the case $\# P > 1$, there is a $z \in P$ with $x < z$ or $z < x$ ($P$ is connected), thus $\tau( \sigma(z)) \in \tau( \sigma( \myouparrow x ))$ or $\tau( \sigma(z)) \in \tau( \sigma( \myodownarrow x ))$ due to the strictness of $\tau \circ \sigma$, in contradiction to $\atsX_3 = \atsX_2 = \emptyset$. Furthermore, the case $\# P = 1$ contradicts $\sigma( y ) \not= \sigma(x)$. Hence, $\atSX$ is the only $\fa \in \alpha_\tau(Q)$ with $\atsX_1 = \fa_1$, $\atsX_2 \subseteq \fa_2$, and $\atsX_3 \subseteq \fa_3$.

\EP

Using the definition of $<_+$, it is easily seen, that the mapping $\phi_\tau$ in \eqref{def_phi_tau} is in fact a strict homomorphism. 

\begin{definition} \label{def_distrsys}
Let $Q' \in \mfP_r$, and let $Q_1, \ldots , Q_L \in \mfP_r$ be a sequence of posets, $L \in \myN_0$. We call a sequence of homomorphisms $\tau_\ell \in \S^o(Q_\ell, Q'), \ell \in \myNk{L}$, a {\em distributor from $Q_1, \ldots , Q_L$ to $Q'$}, iff
\begin{itemize}
\item[(i)] $\tau_\ell$ is distributing for every $\ell \in \myNk{L}$;
\item[(ii)] for every $P \in \mfP^c$, the sets $\T_{\tau_\ell}(P,Q_\ell), \ell \in \myNk{L}$, are pairwise disjoint.
\end{itemize}
\end{definition}

An empty collection of homomorphisms $\tau_\ell \in \S^o(P,Q_\ell)$ with $\ell \in \emptyset$, i.e., with $L = 0$, is always a distributor. For every $Q \in \mfP_r$, there exists a {\em trivial distributor} from $Q$ to $Q$: select $L = 1$, $\tau_1 \in \Aut(Q)$. In the case of $L \geq 2$, it is not required that the posets $Q_\ell$ are pairwise different; for our purpose, the key are the sets $\T_{\tau_\ell}(P,Q_\ell)$, and they are clearly ``sufficiently different'' due to Condition (ii). However, in the case of $L \geq 2$, Condition (ii) enforces $Q' \not= Q_\ell$ for every $\ell \in \myNk{L}$: Assume $L  \geq 2$, $Q_L = Q'$. Then $\tau_L$ is bijective, and Lemma \ref{lemma_bijHom} delivers $\tau_L \in \Aut(Q')$. Let $\ell \in \myNk{L-1}$. For $\sigma \in \S^o(P,Q_\ell)$, $P \in \mfP^c$, we have $\tau_L^{-1} \circ \tau_\ell \circ \sigma \in \S^o(P,Q')$, thus $\tau_\ell \circ \sigma \in \T_{\tau_L}(P,Q_L) \cap \T_{\tau_\ell}(P,Q_\ell)$ in contradiction to (ii).

\begin{corollary} \label{coro_sumqSPQ}
Let $Q' \in \mfP$, let $Q_1, \ldots , Q_L \in \mfP, L \in \myN_0$, and assume that a distributor exists from $Q_1, \ldots , Q_L$ to $Q'$. Then for every $P \in \mfP^c$
\begin{align} \label{ungl_SorPQ}
\sum_{\ell = 1}^L \# \S^o(P, Q_\ell) \quad \leq \quad \# \S^o(P,Q').
\end{align}
\end{corollary}
\BP Let $\tau_\ell \in \S^o(Q_\ell, Q'), \ell \in \myNk{L}$, be a distributor from $Q_1, \ldots , Q_L$ to $Q'$. The sets $\T_{\tau_\ell}(P,Q_\ell)$ are subsets of $\S^o(P,Q')$, and due to (ii) in Definition \ref{def_distrsys}, they are pairwise disjoint. Because $\tau_\ell$ is distributing for every $\ell \in \myNk{L}$, we conclude
\begin{displaymath}
\sum_{\ell=1}^L \# \S^o(P,Q_\ell) \; = \;
\sum_{\ell=1}^L \# \T_{\tau_{\ell}}(P,Q_\ell) \; \leq \;
\# \S^o(P,Q').
\end{displaymath}

\EP

In the conditions of the following theorem, the key are the distributors in (v). The mappings $\lambda_j$ in (ii) select the elements of $\mfE^c(R)$ belonging to the distributor with index $j$, and Condition (iii) ensures that every element of $\mfE^c(R)$ belongs to at least one distributor. The integers $r(i)$ and $q_j(i)$ in (iii) and (iv) count, in how many distributors an element of $\mfE^c( R )$ appears and how often it is repeated in the distributor with index $j$.

\begin{theorem} \label{theo_Transport}
Let $R, S \in \mfP$.
\begin{itemize}
\item[(i)] Let the pairwise different posets $Q_1, \ldots , Q_I \in \mfP_r^c$ form the set $\mfE^c(R)$, and let the pairwise different posets $Q'_1, \ldots , Q'_J \in \mfP_r^c$ form the set $\mfE^c(S)$.
\item[(ii)] Let $\nu_j \in \myN_0$ for every $j \in \myNk{J}$, and define $J^* \equiv \mysetdescr{ j \in \myNk{J}}{ \nu_j \geq 1}$.
For every $j \in J^*$, let $\lambda_j : \myNk{\nu_j} \rightarrow \myNk{I}$ be a mapping.
\item[(iii)] Assume $\myNk{I} = \cup_{j \in J^*} \, \lambda_j \! \left( \myNk{\nu_j} \right)$, and define $r(i) \equiv \# \mysetdescr{ j \in J^* }{ i \in \lambda_j \! \left( \myNk{\nu_j} \right) }$ for every $i \in \myNk{I}$.
\item[(iv)] For every $j \in J^*, i \in \lambda_j \! \left( \myNk{\nu_j} \right)$, define $q_j(i) \equiv \# \myurbild{\lambda_j}(i)$.
\item[(v)] Assume, that a distributor from $Q_{\lambda_j(k)}, k \in \myNk{\nu_j}$, to $Q'_j$ exists for all $j \in J^*$.
\end{itemize}
Then
\begin{align} \label{max_ungl}
\max \mysetdescr{ \frac{ \# \Emb( Q_i, R ) }{ q_j(i) \cdot r(i) \cdot \# \Aut( Q_i ) } }{ i \in \lambda_j \! \left( \myNk{\nu_j} \right) } \quad & \leq \quad \frac{ \# \Emb( Q'_j, S ) }{ \# \Aut( Q'_j )}
\end{align}
for all $j \in J^*$ implies $R \sqsubseteq_G S$.

\end{theorem}
\BP $R \sqsubseteq_G S$ follows with Theorem 1, because for all $P \in \mfP^c$
\begin{align*}
\# \S(P,R) 
& \stackrel{\eqref{eq_summe_bQPT_conn}}{=} \quad \;
\sum_{i=1}^I \# \I_{Q_i}(P,R) \\
& \stackrel{(iii)}{=} \quad \sum_{j \in J^*} \sum_{i \in \lambda_j \! \left( \myNk{\nu_j} \right)} \frac{1}{ r(i ) } \cdot \# \I_{Q_i}(P,R) \\
& \stackrel{(iv)}{=} \quad \;
\sum_{j \in J^*} \sum_{k=1}^{\nu_j} 
\frac{1}{ q_j(\lambda_j(k))\cdot r(\lambda_j(k)) } \cdot \# \I_{Q_{\lambda_j(k)}}(P,R) \\
& \stackrel{\eqref{eq_formel_IQPT_2}}{=} \quad
\sum_{j \in J^*} \sum_{k=1}^{\nu_j} \frac{ \# \S^o(P,Q_{\lambda_j(k)}) }{ q_j(\lambda_j(k)) \cdot r(\lambda_j(k)) } \cdot \frac{ \# \Emb( Q_{\lambda_j(k)}, R )}{ \# \Aut(Q_{\lambda_j(k)}) } \\
& \stackrel{\eqref{max_ungl}}{\leq} \quad 
\sum_{j \in J^*} \left( \sum_{k=1}^{\nu_j} \# \S^o(P,Q_{\lambda_j(k)}) \right) \cdot \frac{ \# \Emb( Q'_j, S )}{ \# \Aut(Q'_j) } \\
& \stackrel{\eqref{ungl_SorPQ}}{\leq} \quad 
\sum_{j \in J^*}  \# \S^o(P, Q'_j) \cdot \frac{ \# \Emb( Q'_j, S )}{ \# \Aut(Q'_j) } \\
& \stackrel{\eqref{eq_formel_IQPT_2}}{=} \quad \;
\sum_{j \in J^*}  \# \I_{Q'_j}(P, S)
\quad \stackrel{\eqref{eq_summe_bQPT_conn}}{\leq} \quad
\# \S(P,S).
\end{align*}

\EP

Corollary \ref{coro_applic_sum_of_soe} is a special case of this theorem: there we have $J = I + K$, $K \in \myN_0$, $Q'_j = Q_j$, $\nu_j = 1$, $\lambda_j = id_{\myNk{1}}$ for all $j \in \myNk{I}$, $\nu_j = 0$ for all $j \in \myNk{J} \setminus \myNk{I}$, and all distributors are trivial.


\section{Application} \label{sec_application}

In Section \ref{subsec_proving_RGS}, we show $R \sqsubseteq_G S$ for two pairs of posets $R$ and $S$ by means of Theorem \ref{theo_Transport}. The reader will realize that both proofs are simple and follow a schematic pattern. In Section \ref{subsec_constr_Gschemes}, we present in Theorem \ref{theo_PQ_AT} a general method of how to construct a poset $T$ for a poset $P + Q$, for which $ P + Q \sqsubseteq_G P \vert_A + T$, where $A$ is a convex subset of $P$. Compared with the difficult proofs for $R \sqsubseteq S$ and $R \sqsubseteq_G S$ in \cite{aCampo_2018,aCampo_toappear_1}, the Theorems \ref{theo_Transport} and \ref{theo_PQ_AT} thus represent a step forward in the development of the efficient toolbox for checking $R \sqsubseteq_G S$ mentioned at the end of Section \ref{sec_introduction}. However, there are cases of $R \sqsubseteq_G S$ in which these two theorems do not help, e.g., Construction (C7) in \cite{aCampo_toappear_1}.

\subsection{Proving $R \sqsubseteq_G S$ by means of Theorem \ref{theo_Transport}} \label{subsec_proving_RGS}

\begin{figure}
\begin{center}
\includegraphics[trim = 70 700 200 70, clip]{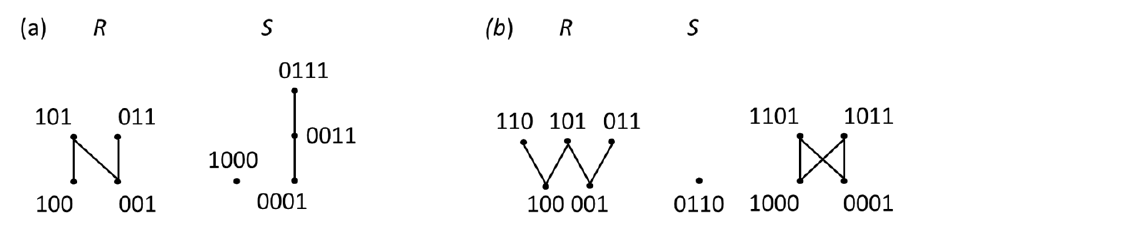}
\caption{\label{fig_NA1C3_WA1N2} Two pairs of posets $R$ and $S$ with $R \sqsubseteq_G S$.}
\end{center}
\end{figure}

Let $R = N$ and $S = A_1 + C_3$, as shown in Figure \ref{fig_NA1C3_WA1N2}(a). By means of Theorem \ref{theo_Transport}, we want to show $R \sqsubseteq_G S$. Looking at Figure \ref{fig_MatrixBspl}, we realize that for all $Q \in \mfE^c(R) \cap \mfE^c(S) = \{ A_1, C_2 \}$, we have $\# \Emb(Q, R) = \# \Emb(Q, S)$. The respective trivial distributors thus fulfill \eqref{max_ungl}. It remains to construct a suitable distributor from $ V_3, \Lambda_3, N \in \mfE^c(R)$ to $C_3 \in \mfE^c(S)$.

\begin{figure}
\begin{center}
\includegraphics[trim = 70 545 200 70, clip]{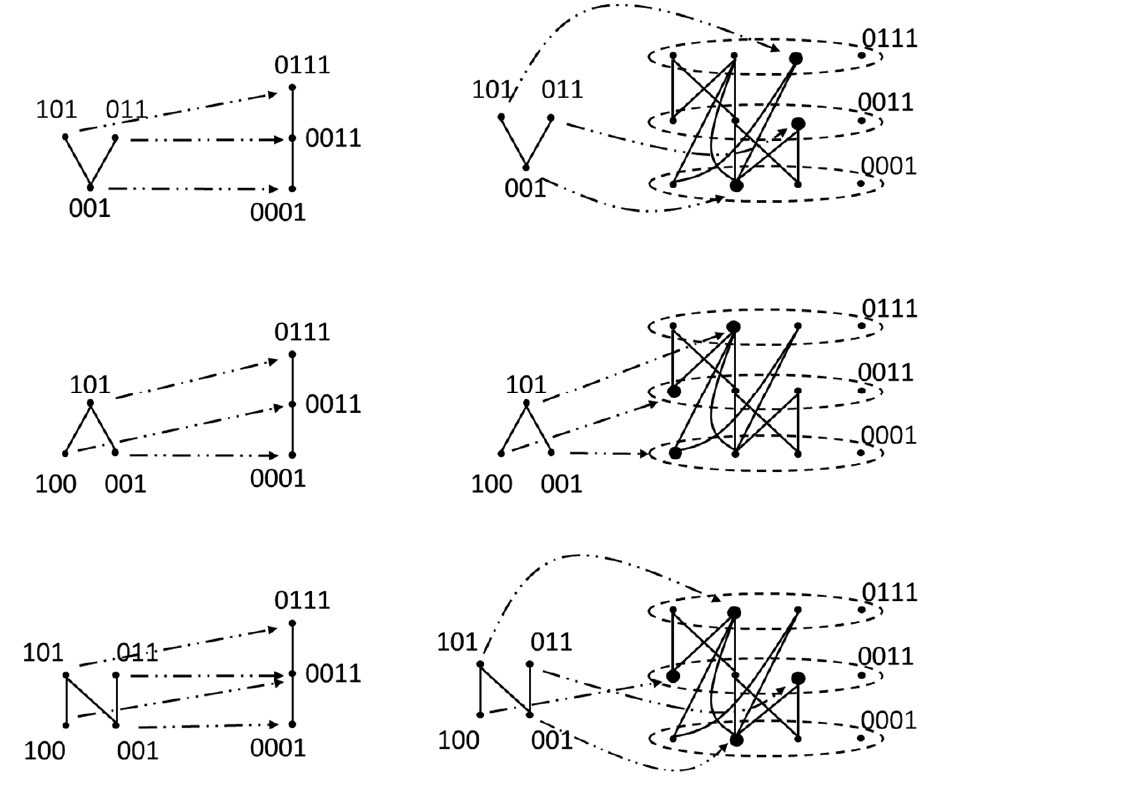}
\caption{\label{fig_N_C3} Top: $\tau_1 \in \S^o( V_3, C_3)$ and $\alpha_{\tau_1} \in \S( V_3, \E(C_3) )$. 
Middle: $\tau_2 \in \S^o( \Lambda_3, C_3 )$ and $\alpha_{\tau_2} \in \S( \Lambda_3, \E(C_3) )$. 
Bottom: $\tau_3 \in \S^o( N, C_3 )$ and $\alpha_{\tau_3} \in \S( N, \E(C_3) )$. }
\end{center}
\end{figure}

In Figure \ref{fig_N_C3}, homomorphisms $\tau_1 \in \S^o( V_3, C_3), \tau_2 \in \S^o( \Lambda_3, C_3 )$, and $\tau_3 \in \S^o( N, C_3 )$ are shown together with their associated homomorphisms $\alpha_{\tau_j}$, $j \in \myNk{3}$. For every $j \in \myNk{3}$, $\alpha_ {\tau_j}$ is one-to-one and \eqref{trennbedingung_2} is fulfilled. According to Lemma \ref{lemma_trennbedingung}(b), the homomorphism $\tau_j$ is thus distributing and fulfills the conditions in Lemma \ref{lemma_trennbedingung}(a). In order to show Condition (ii) in Definition \ref{def_distrsys}, we prove that, among the non-isolated points of $\E(C_3)$, the two outer ones in $\E( C_3; 0011)$ in Figure \ref{fig_N_C3} are characteristic in the following sense: for all $P \in \mfP^c$ and all $\sigma_1 \in \S^o(P,V_3)$, $\sigma_2 \in \S^o(P,\Lambda_3)$, $\sigma_3 \in \S^o(P,N)$
\begin{align*}
\left( 0011, \emptyset, \{ 0111 \} \right) & \not\in \myaxix{\tau_1 \circ \sigma_1}{P} \ni \left( 0011, \{ 0001 \}, \emptyset \right), \\
\left( 0011, \emptyset, \{ 0111 \} \right) & \in \myaxix{\tau_2 \circ \sigma_2}{P} \not\ni \left( 0011, \{ 0001 \}, \emptyset \right), \\
\left( 0011, \emptyset, \{ 0111 \} \right) & \in \myaxix{\tau_3 \circ \sigma_3}{P} \ni \left( 0011, \{ 0001 \}, \emptyset \right).
\end{align*}
This implies that the sets $\T_{\tau_1}(P,V_3), \T_{\tau_2}(P,\Lambda_3)$, and $\T_{\tau_3}(P,N)$ are pairwise disjoint subsets of $\S^o(P,C_3)$ for every $P \in \mfP^c$.

Let $P \in \mfP^c$, $\sigma_1 \in \S^o(P,V_3)$. Because $\tau_1$ fulfills the conditions of Lemma \ref{lemma_trennbedingung}(a), there is no $x \in P$ with $\myaxix{\tau_1 \circ \sigma_1}{x} = \left( 0011, \emptyset, \{ 0111 \} \right)$ (look at Figure \ref{fig_N_C3}). There exists an $x \in P$ with $\sigma_1(x) = 011$. Because $011$ is maximal in $V_3$ and $\sigma_1$ and $\tau_1$ are strict, there is no $z \in P$ with $x < z$, hence $\myaxix{\tau_1 \circ \sigma_1}{x}_3 = \emptyset$. Because $P$ is connected with $\# P \geq 3$, $x$ is not an isolated point There must thus exist a $y \in P$ with $y < x$, thus $001 = \sigma_1(y)$, and we conclude $0001 \in \myaxix{\tau_1 \circ \sigma_1}{x}_2$, hence $\myaxix{\tau_1 \circ \sigma_1}{x}_2 = \{ 0001 \}$.

The proofs of $\left( 0011, \emptyset, \{ 0111 \} \right) \in \myaxix{\tau_2 \circ \sigma_2}{P} \not\ni \left( 0011, \{ 0001 \}, \emptyset \right)$ for all $\sigma_2 \in \S^o(P,\Lambda_3)$ and of $\left( 0011, \emptyset, \{ 0111 \} \right) \in \myaxix{\tau_3 \circ \sigma_3}{P} \ni \left( 0011, \{ 0001 \}, \emptyset \right)$ for all $\sigma_3 \in \S^o(P,N)$ follow exactly the scheme of the proof for $\myaxix{\tau_1 \circ \sigma_1}{x}$.

Each of the posets $Q \in \{ V_3, \Lambda_3, N \}$ occurs exactly once in the sequence $V_3, \Lambda_3, N$, and for each of them we have $\frac{ \# Emb(Q,R) }{ \# \Aut(Q) } = 1 = \frac{ \# Emb(C_3,S) }{ \# \Aut(C_3) }$. Now Theorem \ref{theo_Transport} delivers $ N \sqsubseteq_G A_1 + C_3$.

\begin{figure}[h]
\begin{center}
\includegraphics[trim = 70 480 200 70, clip]{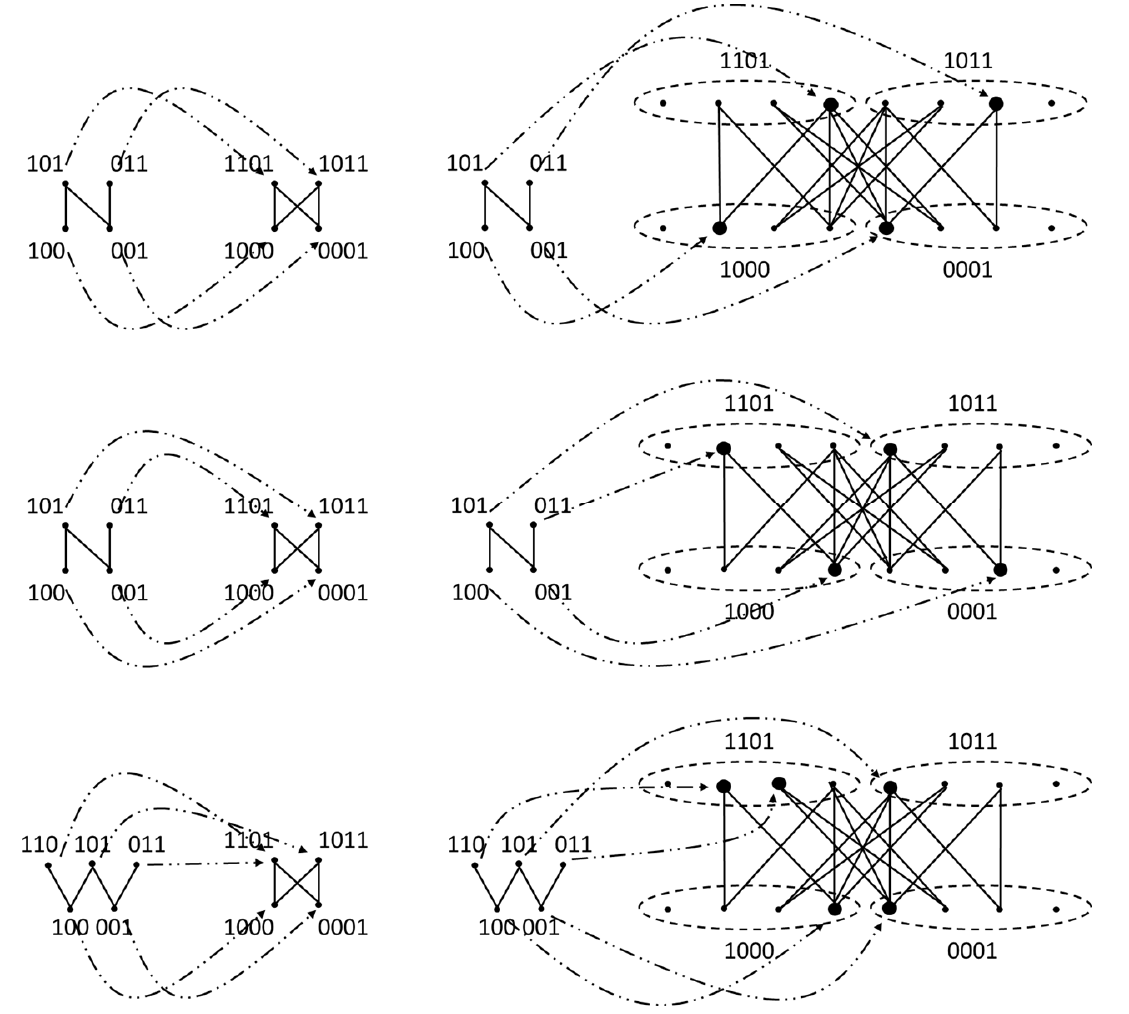}
\caption{\label{fig_W_N2} Top: $\tau_1 \in \S^o( N, N^{(2)})$ and $\alpha_{\tau_1} \in \S( N, \E(N^{(2)}) )$. 
Middle: $\tau_2 \in \S^o( N, N^{(2)} )$ and $\alpha_{\tau_2} \in \S( N, \E(N^{(2)}) )$. 
Bottom: $\tau_3 \in \S^o( W, N^{(2)} )$ and $\alpha_{\tau_3} \in \S( W, \E(N^{(2)}) )$. }
\end{center}
\end{figure}

The second example concerns the posets $R = W$ and $S = A_1 + N^{(2)}$ shown in Figure \ref{fig_NA1C3_WA1N2}(b). According to Figure \ref{fig_MatrixBspl}, we have $\# \Emb(Q, R) \leq \# \Emb(Q, S)$ for all $Q \in \mfE^c(R) \cap \mfE^c(S) = \{ A_1, C_2, V_3, \Lambda_3 \}$. The respective trivial distributors thus fulfill \eqref{max_ungl}. We now construct a distributor of $ N, N, W \in \mfE^c(R)$ to $ N^{(2)} \in \mfE^c(S)$.

In Figure \ref{fig_W_N2}, homomorphisms $\tau_1, \tau_2 \in \S^o( N, N^{(2)})$ and $\tau_3 \in \S^o( W, N^{(2)} )$ are shown together with their respective homomorphism $\alpha_{\tau_j}, j \in \myNk{3}$. All $\alpha_{\tau_j}$ are one-to-one, and all $\tau_j$ fulfill \eqref{trennbedingung_2}. Therefore, all $\tau_j$ are distributing and fulfill the conditions of Lemma \ref{lemma_trennbedingung}(a) according to Lemma \ref{lemma_trennbedingung}(b).

Now it is the innermost points in $\E( N^{(2)}; 1000)$ and in $\E( N^{(2)}; 0001)$ in Figure \ref{fig_W_N2} which make the difference: for every $P \in \mfP^c$, for all $\sigma_1, \sigma_2 \in \S^o(P,N)$, and for all $\sigma_3 \in \S^o(P,W)$
\begin{align*}
\left( 1000, \emptyset, \{ 1101, 1011 \} \right) & \not\in \myaxix{\tau_1 \circ \sigma_1}{P} \ni \left( 0001, \emptyset, \{ 1101, 1011 \} \right), \\
\left( 1000, \emptyset, \{ 1101, 1011 \} \right) & \in \myaxix{\tau_2 \circ \sigma_2}{P} \not\ni \left( 0001, \emptyset, \{ 1101, 1011 \} \right), \\
\left( 1000, \emptyset, \{ 1101, 1011 \} \right) & \in \myaxix{\tau_3 \circ \sigma_3}{P} \ni \left( 0001, \emptyset, \{ 1101, 1011 \} \right), \\
\end{align*}

Let $P \in \mfP^c$ and $\sigma_1 \in \S^o(P,N)$. Because $\tau_1$ fullfills the conditions of Lemma \ref{lemma_trennbedingung}(a), there is no $x \in P$ with $\myaxix{\tau_1 \circ \sigma_1}{x} = \left( 1000, \emptyset, \{ 1101, 1011 \} \right)$ (look at Figure \ref{fig_W_N2}). There exists an $x \in \myurbild{\sigma_1}(101)$ and a $y  \in \myurbild{\sigma_1}(011)$, and there is a zigzag line $z_0, \ldots , z_L$ connecting $x$ and $y$. We conclude that there is an $\ell \in \myNk{L-1}$ with $\sigma_1(z_{\ell-1}) = 101$, $\sigma_1(z_\ell) = 001$, and $\sigma_1(z_{\ell+1}) = 011$, hence $\myaxix{\tau_1 \circ \sigma_1}{z_\ell}_3 = \{ 1101, 1011 \}$. Furthermore, $001$ is a minimal point in $N$ and $\sigma_1$ is strict. There is thus no $a \in P$ with $a < z_\ell$, and $\myaxix{\tau_1 \circ \sigma_1}{z_\ell}_2 = \emptyset$ is shown.

The proof of the claim about $\myaxix{\tau_2 \circ \sigma_2}{P}$ runs exactly the same way, and the claim about $\myaxix{\tau_3 \circ \sigma_3}{P}$ is proven by applying the pattern of the second part of the proof for $\myaxix{\tau_1 \circ \sigma_1}{x}$ twice.

Because of $\frac{ \# \Emb( N, W ) }{ 2 \cdot \# \Aut(N) } = 1$, $\frac{ \# \Emb( W, W ) }{ \# \Aut(W) } = 1$, and $\frac{ \# \Emb( N^{(2)}, N^{(2)} ) }{ \# \Aut(N^{(2)}) } = 1$, Theorem \ref{theo_Transport} yields $W \sqsubseteq_G A_1 + N^{(2)}$.

\subsection{A construction method based on Corollary \ref{coro_applic_sum_of_soe}} \label{subsec_constr_Gschemes}

\begin{figure}[h]
\begin{center}
\includegraphics[trim = 70 730 200 70, clip]{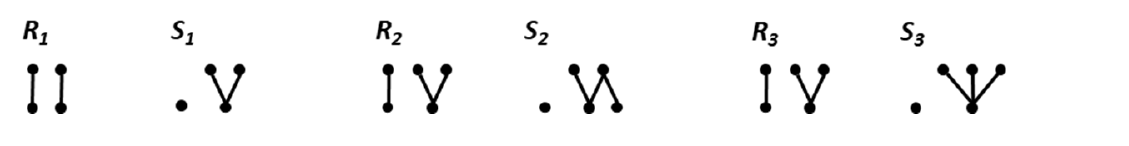}
\caption{\label{fig_ConstrBspl} Three pairs of posets $R_i, S_i$ with $R_i \sqsubseteq_G S_i$.}
\end{center}
\end{figure}

\begin{figure}
\begin{center}
\includegraphics[trim = 70 700 200 70, clip]{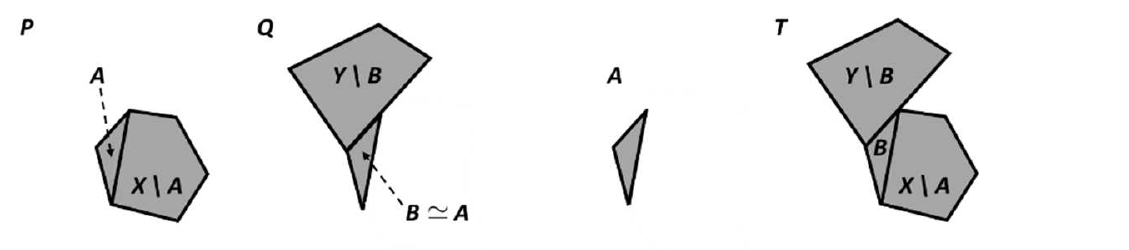}
\caption{\label{fig_ConstrMeth} Illustration of a method for the construction of posets $R = P + Q$ and $S = P \vert_A + T$ with $R \sqsubseteq_G S$.}
\end{center}
\end{figure}

In Figure \ref{fig_ConstrBspl}, three pairs of posets $R, S$ are shown with $R \sqsubseteq_G S$ according to \cite{aCampo_toappear_1}. In all pairs, the poset $R$ has the form $P + Q$ and the poset $S$ has the form $A + T$ where $A$ consists of a single point only. In this section we develop a construction method covering such pairs of posets; the key element in showing $R \sqsubseteq_G S$ will be Corollary \ref{coro_applic_sum_of_soe}.

The construction method works as shown in Figure \ref{fig_ConstrMeth}. Let $P, Q$ be disjoint posets, and let $A$ be convex in $P$. We select a convex subset $B$ of $Q$ with $Q \vert_B \simeq P \vert_A$, and implement $P \vert_{X \setminus A}$ into $Q$ by connecting it with $B$ as it was connected with $A$ in $P$. Denoting with $T$ the poset resulting from this combination of $P \vert_{X \setminus A}$ and $Q$, we show $P + Q \sqsubseteq_G P \vert_A + T$ in Theorem \ref{theo_PQ_AT}. (In Figure \ref{fig_ConstrBspl}, all $R$ have the form $C_2 + Q$, and it is in all three cases one point of the $C_2$-component of $R$ which is implemented into $Q$ whereas the remaining point forms the subset $A$.) The approach is formalized as follows:

\begin{definition} \label{def_leq_T}
In this section, we agree on the following:
\begin{itemize}
\item $P = ( X, \leq_P )$ and $Q = ( Y, \leq_Q )$ are disjoint finite posets.
\item $A \subseteq X$ is convex in $P$, $B \subseteq Y$ is convex in $Q$, and $ \beta : P \vert_A \rightarrow Q \vert_B$ is an isomorphism.
\item $W \equiv X \setminus A$ and $O \equiv P \vert_W $; the partial order relation of $O$ is denoted by $\leq_O$.
\item We define on $V \equiv W \cup Y$ a binary relation $\leq_T$ by setting
\end{itemize}
\begin{align*}
\leq_T & \equiv \; \; \leq_O \; \cup \; \leq_Q \; \cup \; \leq_d \; \cup \;  \leq_u \\
\mytext{where} \; \leq_d & \equiv \bigcup_{w \in W}
\mysetdescr{ ( y, w ) \in Y \mytimes W }{ \exists \, a \in A : \, (a,w) \inA \leq_P \mytext{and} y \in  \downarrow_Q \beta(a) }, \\
\leq_u & \equiv \bigcup_{w \in W}
\mysetdescr{ ( w, y ) \in W \mytimes Y}{ \exists  \, a \in A : \, (w,a) \inA \leq_P \mytext{and} y \in \uparrow_Q \beta(a) }.
\end{align*}
\end{definition}

The relations $\leq_O$, $\leq_Q$, $\leq_d$, and $\leq_u$ are pairwise disjoint. We will see in Lemma \ref{lemma_leqT_is_partOrd} that $T \equiv (V, \leq_T)$ is a partial order. The subsets $\leq_O$ and $\leq_Q$ transfer the partial orders of $O$ and $Q$ to $T$, and the addition of $\leq_d$ and $\leq_u$ fixes $O$ on $B$ and makes $\leq_T$ transitive. In the most simple case, $A$ and $B$ are antichains in $P$ and $Q$, respectively, with $ \# A = \# B$, as in Figure \ref{fig_ConstrBspl}.

The proof of the next lemma requires some technical work in which the following equivalences are tacitly applied: For every $(v,w) \inA \leq_T$
\begin{align*}
v \in W & \Leftrightarrow (v,w) \inA \leq_O \; \cup \leq_u, \\
v \in Y & \Leftrightarrow (v,w) \inA \leq_Q \; \cup \leq_d, \\
w \in W & \Leftrightarrow (v,w) \inA \leq_O \; \cup \leq_d, \\
w \in Y & \Leftrightarrow (v,w) \inA \leq_Q \; \cup \leq_u.
\end{align*}

\begin{lemma} \label{lemma_leqT_is_partOrd}
$\leq_T$ is a partial order on $V$. The pair $T = ( V, \leq_T)$ is thus a poset.
\end{lemma}
\BP {\em Reflectivity:} Trivial due to $\leq_O \cup \leq_Q \subseteqA \leq_T$.

{\em Antisymmetry:} Let $v, w \in V$ with $(v,w), (w,v) \inA \leq_T$.
\begin{itemize}
\item $(v,w) \inA \leq_O$: Then $v, w \in W$, thus $(w,v) \inA \leq_O$, and $v = w$ follows.
\item $(v,w) \inA \leq_Q$: Proven as the case $(v,w) \inA \leq_O$.
\item $(v,w) \inA \leq_d$: Then $v \in Y, w \in W$, hence $(w,v) \inA \leq_u$. There exist thus $a, b \in A$ with
\begin{align*}
(a,w) \inA \leq_P & \mytext{and} v \inA \downarrow_Q \beta(a) 
\quad \mytext{due to} (v,w) \inA \leq_d,\\
(w,b) \inA \leq_P & \mytext{and} v \in \uparrow_Q \beta(b)
\quad \; \mytext{due to} (w,v) \inA \leq_u.
\end{align*}
But then $a \leq_P w \leq_P b$, thus $w \in A$, because $A$ is convex. But this contradicts $w \in W$.
\item $(v,w) \inA \leq_u$: Proven as the case $(v,w) \inA \leq_d$.
\end{itemize}
{\em Transitivity:} Let $(v,w), (w,z) \inA \leq_T$.
\begin{itemize}
\item $(v,w) \inA \leq_O$: Then we have $w \in W$, thus $(w,z) \inA \leq_O \cup \leq_u$.

In the case of $(w,z) \inA \leq_O$, we have $(v,z) \inA \leq_O$.

Let $(w,z) \inA \leq_u$. There is an $a \in A$ with $(w,a) \inA \leq_P$ and $z \inA \uparrow_Q \beta(a)$. According to our assumption, we have $(v,w) \inA \leq_O \subseteqA \leq_P$, which yields $(v,a) \inA \leq_P$ with $z \inA \uparrow_Q \beta(a)$. We conclude $(v,z) \inA \leq_u$.

\item $(v,w) \inA \leq_Q$: Proven as the case $(v,w) \inA \leq_O$.

\item $(v,w) \inA \leq_d$: Then we have $w \in W$, thus $(w,z) \inA \leq_O \cup \leq_u$. Due to $(v,w) \inA \leq_d$, there exists an $a \in A$ with $(a,w) \inA \leq_P$ and $v \in \downarrow_Q \beta(a)$.

In the case $(w,z) \inA \leq_O \subseteqA \leq_P$, we get $(a,z) \inA \leq_P$, hence $(v,z) \inA \leq_d$.

In the case $(w,z) \inA \leq_u$, there exists a point $b \in A$ with $(w,b) \inA \leq_P$ and $z \in \uparrow_Q \beta(b)$. With $(a,w) \inA \leq_P$, we conclude $(a,b) \inA \leq_P$, thus $(\beta(a), \beta(b))$ $ \inA \leq_Q$. Now we get $v \leq_Q \beta(a) \leq_Q \beta(b) \leq_Q z$, thus $(v,z) \in \leq_Q$.

\item $(v,w) \inA \leq_u$: Proven as the case $(v,w) \inA \leq_d$.

\end{itemize}

\EP

\begin{theorem} \label{theo_PQ_AT}
With the construction in Definition \ref{def_leq_T}, $P + Q \sqsubseteq_G P \vert_A + T$.
\end{theorem}
\BP With $P' \equiv T \vert_{W \cup B}$ and $Q' \equiv T \vert_Y$, we have $P \simeq P'$ and $Q = Q'$. (For $P \simeq P'$, observe that for all $a \in A, w \in W$, $(a,w) \inA \leq_Q$ implies $(\beta(a),w) \inA \leq_d$ and $(w,a) \inA \leq_Q$ implies $(w,\beta(a)) \inA \leq_u$. Furthermore, $\beta : P \vert_A \rightarrow Q \vert_B$ is an isomorphism.)

In the proof, we will use Corollary \ref{coro_applic_sum_of_soe}. With $A' \equiv P \vert_A$, we have thus to show $\mfE^c(P+Q) \subseteq \mfE^c(A'+T)$ and 
$ \# \Emb(E,P+Q) \leq \# \Emb(E,A'+T)$ for every $E\in \mfE^c(P+Q)$.

Because $\mfE^c(P + Q)$ contains connected posets only, we get
\begin{align*}
\mfE^c(P + Q) & = \mfE^c(P) \cup \mfE^c(Q) = \mfE^c(P') \cup \mfE^c(Q')
\subseteq \mfE^c(T) \subseteq \mfE^c(A' + T).
\end{align*}

Let $E \in \mfE^c(P + Q)$. Assume that for $\xi \in \Emb(E, P' )$ and $\zeta \in \Emb(E, Q' )$, we have $\xi(x) = \zeta(x)$ for all $x \in E$, i.e., $\xi \vert^V = \zeta \vert^V$. Then $\xi(E) = \xi(E) \cap \zeta(E) \subseteq (W \cup B) \cap Y = B$, hence $\beta^{-1} \circ \xi \in \Emb(E,A')$. With the sets 
\begin{align*}
F_1 & \equiv \mysetdescr{ \xi \vert^V }{ \xi \in \Emb(E, P') } \; \subseteq \; \Emb(E,T), \\
F_2 & \equiv \mysetdescr{ \zeta \vert^V }{ \zeta \in \Emb(E, Q') } \; \subseteq \; \Emb(E,T),
\end{align*}
the mapping $F_1 \cap F_2 \rightarrow \Emb(E,A')$, $\xi \mapsto \beta^{-1} \circ \xi$ is thus well-defined and one-to-one. Because $E$ is connected, we conclude
\begin{align*}
\# \Emb( E, A' + T ) & = \# \Emb( E, A' ) + \# \Emb( E, T ) \\
& \geq \# \Emb( E, A' ) + \# \left( F_1 \cup F_2 \right) \\
& \geq \# F_1 + \# F_2 \\
& = \# \Emb(E, P') + \# \Emb(E, Q') \\
& = \# \Emb(E, P+Q).
\end{align*}
Now apply Corollary \ref{coro_applic_sum_of_soe}.

\EP

If $A$ is an antichain in $P$, we achieve an even stronger result. With $A' \equiv P \vert_A$, let $\psi : P \rightarrow A' + T$ be the embedding of $P$ into $T$:
\begin{align*}
\psi(x) & \equiv 
\begin{cases}
x, & \mytext{if} x \in X \setminus A; \\
\beta(x), & \mytext{if} x \in A.
\end{cases}
\end{align*}
We have $\E(P+Q) = \E(P) + \E(Q)$, $\E(A' + T) = \E(A') + \E(T)$, and $\E( A' ) = \mysetdescr{ ( a, \emptyset, \emptyset ) }{ a \in A }$. We define a mapping $\eps : \E(P + Q) \rightarrow \E( A' + T)$ by
\begin{align*}
\eps(\fa) & \equiv 
\begin{cases}
\fa, & \mytext{if} \fa \in \E(Q) \cup \E(A'); \\
\left( \psi( \fa_1 ), \psi( \fa_2 ), \psi( \fa_3 ) \right), & \mytext{if} \fa \in \E(P) \setminus \E(A').
\end{cases}
\end{align*}
Because $A$ is an antichain in $P$, the mapping $\eps$ is a one-to-one homomorphism. A short calculation shows that $\eps$ fulfills all conditions of \cite[Theorem 5]{aCampo_toappear_1}, and we conclude that there exists a strong I-scheme from $P + Q$ to $A' + T$. (A strong I-scheme is a strong G-scheme fulfilling an additional regularity condition; strong I-schemes are treated in \cite{aCampo_toappear_1}.)

\end{document}